\documentclass{article}

\usepackage{amssymb} 
\usepackage{amsmath} 
\usepackage{latexsym} 
\usepackage{theorem} 
\usepackage{psfrag}
\usepackage{epsfig}
\usepackage[latin1]{inputenc}

\def\vect#1{\mbox{\boldmath $#1$}}

\theorembodyfont{\itshape} 

\newtheorem{theorem}{Theorem}[section]
\newtheorem{lemma}[theorem]{Lemma}

\newtheorem{corollary}[theorem]{Corollary}
\newtheorem{proposition}[theorem]{Proposition}

\newtheorem{definition}[theorem]{Definition}

\newenvironment{remark}{{\noindent \bf Remark. }}{}
\newenvironment{proof}{{\par\addvspace{0.1cm}\noindent \bf Proof. }}{\hfill$\Box$\par\medskip}

\newcommand{\spa}{\mathrm{span}}
\def\span#1{\spa\,\langle #1 \rangle}
\def\spanbig#1{\spa\,\big\langle #1 \big\rangle}

\def\w{\wedge}
\def\p{\prime}
\def\La{\Lambda}
\def\Si{\Sigma}

\def\omegasub#1{\mbox{\large $\omega$}\mbox{\small${\mbox{\large ${}$}}_{#1}$}}
\def\spbmapright#1#2{\smash{%
 \mathop{\hbox to 0.8cm{\rightarrowfill}}
  \limits^{#1}_{#2}}}

 \def\ppt#1{{\overset{\mbox{\large .\kern-1pt.}}{#1}}}
 \def\pppt#1{{\overset{\mbox{\large .\kern-1pt.\kern-1pt.}}{#1}}}

 \def\Pt#1{{\overset{\mbox{\Huge .}}{#1}}}
 \def\Ppt#1{{\overset{\mbox{\Huge .\kern-3pt.}}{#1}}}
 \def\Pppt#1{{\overset{\mbox{\Huge .\kern-3pt.\kern-3pt.}}{#1}}}
\def\HH{\mathbb{H}}

\def\RR{\mathbb{R}}
\def\EE{\mathbb{E}}
\def\CC{\mathbb{C}}
\def\ss{{\bf S}}
\def\LL{{\mathcal L}}
\def\SS{\vect{S}}
\def\Ga{\Gamma}

\def\lc{\mbox{{$\mathcal{L}$}\it ight}} 

\usepackage{color}

\def\ts{\tilde{s}}

\numberwithin{equation}{section}

\title{Conformal arc-length as $\frac12$ dimensional length of the set of osculating circles}

\author{R\'emi Langevin and Jun O'Hara\footnote{This work is partly supported by the JSPS (Japan Society for the Promotion of Science) Bilateral Program and Grant-in-Aid for Scientific Research No. 19540096.}\\
{\small Institut de Math\'ematiques de Bourgogne, Universit\'e de Bourgogne}\\
{\small Department of Mathematics, Tokyo Metropolitan University}}

\begin{document}

\maketitle

\begin{abstract} 
The set of osculating circles of a given curve in $\SS^3$ forms a lightlike curve in the set of oriented circles in $\SS^3$. 
We show that its ``${\frac12}$-dimensional measure'' with respect to the pseudo-Riemannian structure of the set of circles is proportional to the conformal arc-length of the original curve, which is a conformally invariant local quantity discovered in the first half of the last century. 
\end{abstract}



\section{Introduction}
The Frenet-Serret formula provides a local expression of a space curve in terms of the arc-length, the curvature, and the torsion. 
It is well-known that a space curve is determined by the curvature and the torsion up to motion of the Euclidean space $\RR^3$, i.e. isometric transformation of $\RR^3$. 

Let us consider local theory of space curves in conformal geometry. 
We remark that the arc-length is not preserved by M\"obius transformations. 
Three conformal invariants have been found using a suitable normal form. 
They are conformal arc-length, conformal curvature, and conformal torsion (the reader is referred to \cite{CSW} for example). 
Just like in the Euclidean case we have: 
\begin{theorem}{\rm (\cite{Fi}, Theorem 7.2)} 
An oriented connected vertex-free curve is determined up to conformal motion by the three conformal invariants, the conformal arc-length, the conformal torsion, and the conformal curvature. 
\end{theorem}

In this article, we study the conformal arc-length.
\begin{definition}\rm 
Let ${C}$ be an oriented curve in $\mathbb R^3$. 
Let $s, \kappa, \tau$ be the arc-length, curvature, and torsion of ${C}$ respectively. 
The {\em conformal arc-length} parameter $\rho$ of ${C}$ is given by 
\begin{eqnarray}\label{inf_conf_arc-length}
d\rho=
\sqrt[4]{{\kappa^{\p}}^2+\kappa^2\tau^2\,}\,ds. 
\end{eqnarray}
It gives a conformally invariant parametrization of a vertex-free curve. 
We call the $1$-form $(c^{-1})^{\ast}d\rho$ on the curve $C$ the {\em conformal arc-length element}, where $c$ is a map from some interval $I$ to $\RR^3$ so that $C=\{c(s)\}$. 
\end{definition}

The conformal arc-length was given in \cite{Li} and the above formula was given in \cite{Ta}. 

\smallskip
In this paper, we give a new interpretation of the conformal arc-length in terms of the set of the osculating circles. 

\smallskip
Let $\gamma$ be a lightlike curve. 
Although its length 
is equal to $0$, 
we can define a non-trivial ``{$L^{\frac12}$-measure}" of $\gamma$ by 
\[L^{\frac12}(\gamma)=\lim_{\max|t_{j+1}-t_j|\to+0}\sum_{i}\sqrt{\Vert\gamma(t_{i+1})-\gamma(t_i)\Vert}\,,\]
and a ``$\frac12$ dimensional length element'' $d\rho_{\!L^{\frac12}(\gamma)}$ by 
$\displaystyle d\rho_{\!L^{\frac12}(\gamma)}=\sqrt[4]{\frac{\big|\big\langle\Ppt{\gamma}, \Ppt{\gamma}\big\rangle\big|}{12}} \,\,dt$ 
so that $\displaystyle L^{\frac12}(\gamma)=\int_{\gamma}d\rho_{\!L^{\frac12}(\gamma)}\,$. 

Let $\mathcal{S}(1,3)$ denote the set of the oriented circles in $\RR^3$ (or $\SS^3$), where we consider lines in $\RR^3$  as circles. 
It has a pseudo-Riemannian structure with index $2$ which is compatible with the M\"obius transformations. 
Let $C$ be a curve in $\RR^3$. 
The set of osculating circles to $C$ forms a {lightlike} curve $\gamma$ in $\mathcal{S}(1,3)$ (Theorem \ref{thm_osc_circles}). 
Our main theorem claims that the $L^{\frac12}$-measure of $\gamma$ 
is equal to a constant times the conformal arc-length of the original curve $C$ (Corollary \ref{thm_osc_circles2}). 

We also study various properties of curves of osculating circles, or in general, lightlike  curves in the set of oriented circles (or spheres). 
{In each case, a geometric counterpart of the $L^{\frac{1}{2}}$ measure  is given.}  

\smallskip
This article is arranged as follows. 
In section \ref{sec_preliminaries} we explain how to realize spheres and Euclidean spaces in Minkowski space. 
In section \ref{sec_how_to} we give the bijection between the set of codimension $1$ oriented spheres and the de Sitter space. 
In section \ref{sec_half_dim_length} we define $\frac12$ dimensional length element and $L^{\frac12}$-measure of a lightlike curve. 
In section \ref{sec_lightlike_curves} we study the set of osculating circles to a curve in $\RR^2$ (or $\SS^2$), which becomes a lightlike curve in $3$-dimensional de Sitter space.  
The section \ref{sec_S(1,3)} is a preparation to the section \ref{sec_osculating_circles}, where we prove our main theorem. 
These two sections can be read independently of the previous section. 
In the last section we relate lightlike curves in the space of circles,  lightlike curves in the space of spheres and the infinitesimal cross ratio. 
\smallskip

In this article we will use the following notations: 
the letter $s$ is used for the arc-length of a curve in $\EE^3$ or $\SS^3$, 
and the derivative with respect to $s$ is denoted by putting $\prime$$\,$. 
The letter $t$ is for a general parameter, and the derivative with respect to $t$ is denoted by putting $\Pt{}$. The letter $\ts$ is for the arc-length of the set of osculating spheres to a curve in $\EE^3$ or $\SS^3$, which is a curve in de Sitter space in the Minkowski space (Section \ref{sec_last_section}). 

\smallskip
The authors thank Gil Solanes and Martin Guest for valuable discussions. 

\section{Spherical and Euclidean models in the Minkowski space}\label{sec_preliminaries}
%
%
Let $n$ be $1,2$, or $3$. 
The {\em Minkowski space} $\mathbb{R}^{n+2}_1$ is $\mathbb{R}^{n+2}$ with indefinite inner product: 
$$\langle \vect x, \vect y \rangle\!=\!-x_0y_0+x_1y_1+\cdots+x_{n+1}y_{n+1}.$$

Define the {\em Lorentz form} by $\LL(\vect v)=\langle\vect v, \vect v\rangle$. 
The {\em norm} of a vector $\vect v$ is given by $\Vert\vect v\Vert=\sqrt{|\LL(\vect v)|\,}$. 
A vector $\vect v$ in $\mathbb{R}^{n+2}_{1}$ is called {\em spacelike} if $\LL(\vect v)>0$, {\em lightlike} if $\LL(\vect v)=0$ and $\vect v\ne\vect 0$, and {\em timelike} if $\LL(\vect v)<0$. 
The set of lightlike vectors and the origin, $\left.\left\{\vect v\in\mathbb{R}^{n+2}_1 \, \right|\, \langle\vect v,\vect v\rangle=0\right\}$, is called the {\em light cone} and shall be denoted by $\lc$. 
The ``pseudo-sphere'', $\{\vect v\in\mathbb{R}^{n+2}_1|\langle\vect v,\vect v\rangle=1\}$, is called the {\em de Sitter space} and shall be denoted by $\La$ or $\La^{n+1}$. 
\begin{figure}[htbp]
\begin{center}
\includegraphics[width=65mm]{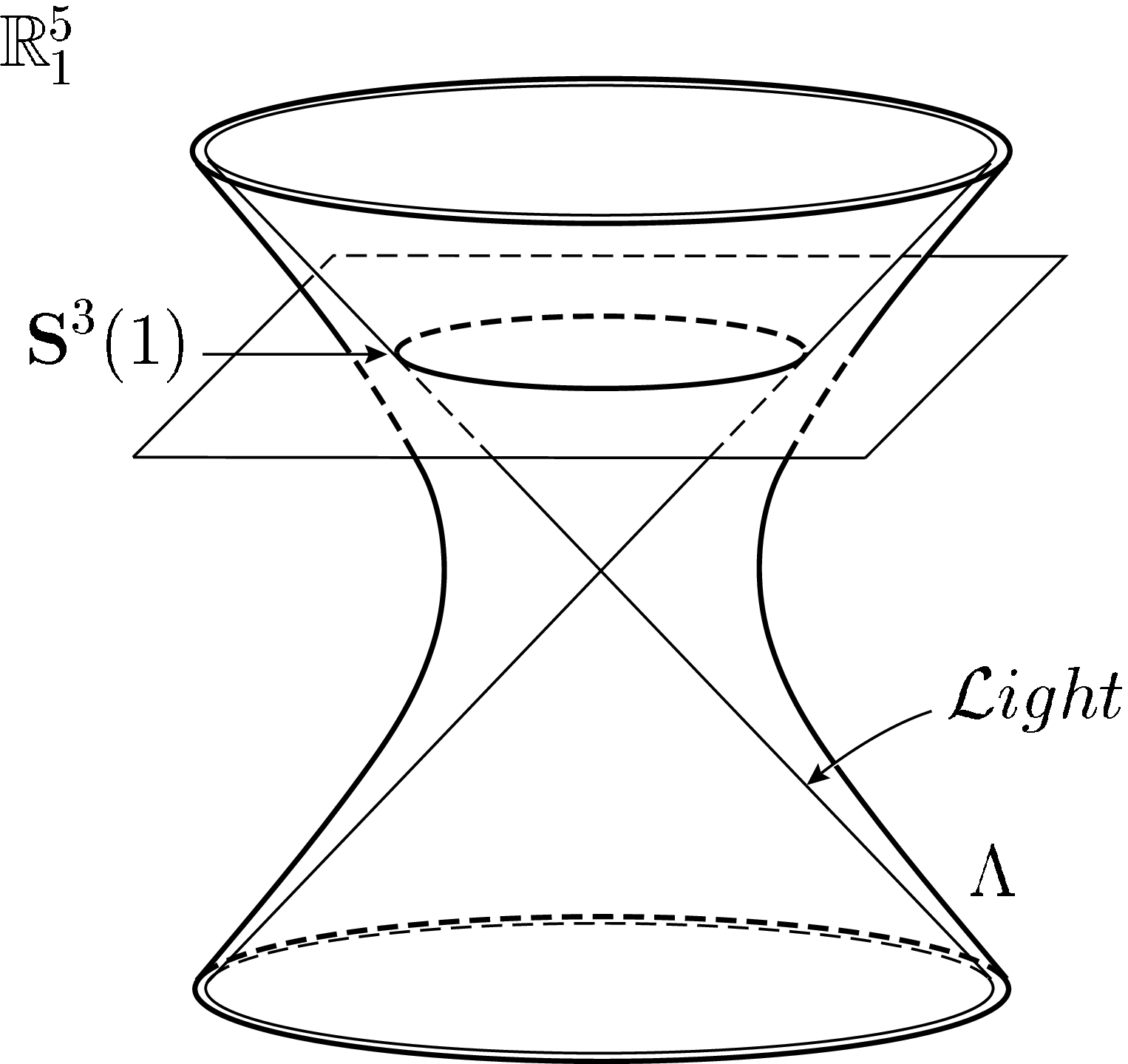}
\caption{Model of $\ss^3$, light cone and de Sitter space \label{quadric-ce3}}
\end{center}
\end{figure}

\begin{figure}[htb]
\begin{center}
\includegraphics[width=.99\linewidth]{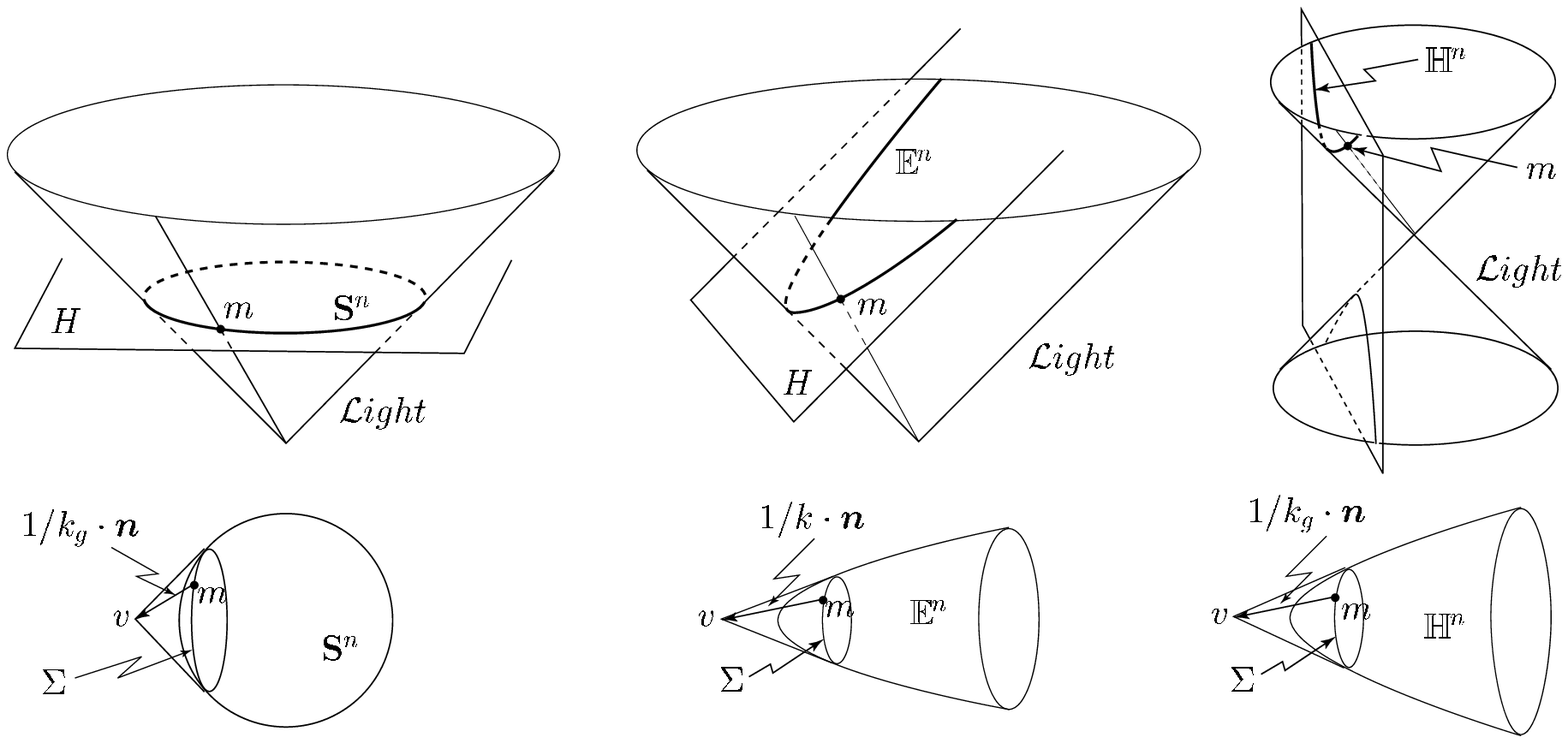}
\caption{Spherical, Euclidean, and hyperbolic models in the Minkowski space $\mathbb R^{n+2}_1$ (up). The geodesic curvature $k_g$, pictures in affine hyperplanes $H$ (down). }
\label{rad_curv_geod}
\end{center}
\end{figure}

Let $W$ be a vector subspace of $\mathbb{R}^{n+2}_{1}$. 
There are three cases which are mutually exclusive. 
Let $\langle\,,\,\rangle|_{W}$ denote the restriction of $\langle \,,\,\rangle$ to $W$. 
\begin{enumerate}
\item The case when $\langle\,,\,\rangle|_{W}$ is non-degenerate. 
This case can be divided into two cases: 
\begin{enumerate}
\item[(1-a)] The case when $\langle\,,\,\rangle|_W$ is indefinite. 
It happens if and only if $W$ intersects the light cone transversely. 
In this case $W$ is said to be {\em timelike}. 
\item [(1-b)] The case when $\langle\,,\,\rangle|_W$ is positive definite. 
It happens if and only if $W$ intersects the light cone only at the origin. 
In this case $W$ is said to be {\em spacelike}. 
\end{enumerate}
\item The case when $\langle\,,\,\rangle|_W$ is degenerate. 
It happens if and only if $W$ is tangent to the light cone. 
In this case $W$ is said to be {\em isotropic}. 
\end{enumerate}
The sphere $\SS^n$, the Euclidean space $\EE^n$, and the hyperbolic space $\HH^n$ can be realized in $\RR_1^{n+2}$ as affine sections of the light cone, i.e. the intersection of an affine $(n+1)$-space $H$ and the light cone (Figure \ref{rad_curv_geod} up). 
We call them spherical, Euclidean, and hyperbolic models respectively. 
Their metrics are induced from the Lorentz form on $\RR_1^{n+2}$ (\cite{G}). 

(1) When the affine space $H$ is tangent to the hyperboloid $F_+=\{\vect x\,|\,\LL(\vect x) = -1, x_0>0\}$, the intersection $\lc \cap H$ is a sphere $\SS^n$ with constant curvature $1$. 

When the tangent point is $(1,0,\cdots,0)$, i.e. when $H=\{x_0=1\}$ we will denote it by $\SS^n(1)$. 
The $n$-sphere $\SS^3$ can be identified with the set of lines through the origin in the light cone. 

(2) When the affine space $H$ is parallel to an isotropic subspace and does not contain the origin, the intersection $\lc \cap H$ is an Euclidean space $\EE^n$. 
For example, take two lightlike vectors $\vect n_1$ and $\vect n_2$ given by 
\[\vect n_1=(1,1,0,\cdots,0), \,\vect n_2=(1,-1,0,\cdots,0).\]
Put 
\[H=\vect n_2+\left(\textsl{Span}\langle\vect n_1\rangle\right)^{\perp} \>\>\mbox{\rm and }\>\>\> \EE^n_0=\lc\cap H.\]
Then the intersection can be explicitly expressed as 
\begin{equation}\label{def_E^3_0}
\EE^n_0=\left\{\left.\left(1+\frac{\vect x\cdot \vect x}4,\,-1+\frac{\vect x\cdot \vect x}4,\, \vect x\right)\,\right|\,\vect x\in \RR^n\right\}\,,
\end{equation}
where $\cdot$ denotes the standard inner product. 
There is an isometry obtained by the projection in the direction of $\vect n_1$ from $\EE^n_0$ to $\{(1,-1,\vect x)\,|\,\vect x\in\RR^n\}\cong\RR^n$. 

The lightlike lines in the light cone gives the bijection between $\SS^n(1)\setminus\{\vect n_1\}$ and $\EE^n_0$. 
It is exactly same as the stereographic projection from the north pole $N$ from $\SS^n\setminus\{N\}$ to $\RR^n$ which is tangent to $\SS^n$ at the south pole through the identifications $\SS^n\cong\SS^n(1)$ and $\RR^n\cong\EE^n_0$. 
\section{de Sitter space as the space of codimension $1$ spheres}\label{sec_how_to}
%
Let $\mathcal{S}(n-1,n)$ be the set of oriented $(n-1)$-spheres $\Si$ in $\EE^n$ (or $\SS^n$). 
Then there is a bijection from $\mathcal{S}(n-1,n)$ to the de Sitter space $\Lambda^{n+1}$. 
Let us express this bijection $\varphi$ in two ways. 
In this section we consider $\EE^n$ (or $\SS^n$) as the intersection of the light cone and an affine hyperplane $H$ in the Minkowski space $\mathbb R^{n+2}_1$. 
\subsection{Using pseudo-orthogonality} 
\begin{figure}[htbp]
\begin{center}
\includegraphics[width=.4\linewidth]{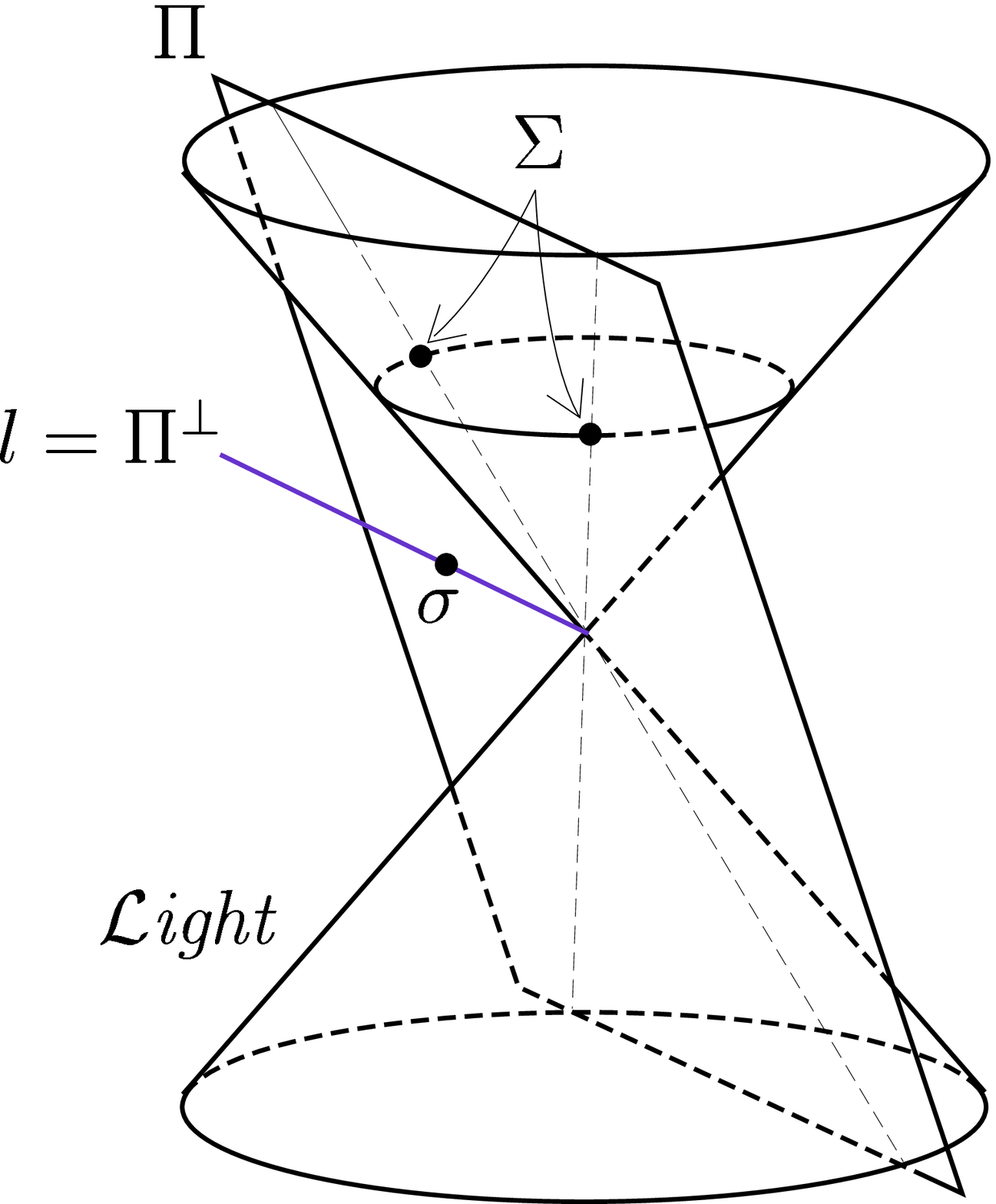}
\caption{The bijection between $\mathcal{S}(n-1,n)$ and $\Lambda^{n+1}$\label{Lambda-c-e}}
\end{center}
\end{figure}

Let $\Si$ be an oriented $(n-1)$-sphere in $\EE^n$ (or $\SS^n$). 
Then $\Si$ can be obtained in $H$ as the intersection of $\EE^n$ (or $\SS^n$) and an affine hyperplane $W$ of $H$. 
By taking a cone from the origin of $\mathbb{R}^{n+2}_1$, $\Si$ can be realized in $\mathbb{R}^{n+2}_1$ as the intersecton of the light cone and an oriented codimension $1$ vector subspace of $\mathbb{R}^{n+2}_1$, $\Pi$ (Figure \ref{Lambda-c-e}). 
Let $\sigma\in\La^{n+1}$ be the endpoint of the positive unit normal vector to $\Pi$. 
\begin{definition} \rm 
As above, the bijection $\varphi: \mathcal{S}(n-1,n)\to\La^{n+1}$ is given by assigning $\sigma$ to $\Sigma$. 
\end{definition}

Since the pseudo-orthogonality is preserved by the action of the Lorentz group $O(4,1)$, this bijection is compatible with the action of $O(4,1)$, i.e. 
\begin{equation}\label{conf_equivar_varphi}
\varphi\,(A\cdot\Sigma)=A\,\varphi(\Sigma)\hspace{0.5cm}(A\in O(4,1)). 
\end{equation}

\subsection{Using geodesic curvature and a normal vector}
Let $\Si$ be a codimension $1$ sphere in $\SS^n$ or in $\EE^n$. 
The {\em geodesic curvature} of a sphere $\Si$ means the geodesic curvature of any geodesic curve on it. 
Let it be denoted by $k_g$. 
Let $m$ be a point in $\Si$ and $\vect n$ the unit normal vector to $\Si$ at $m$. 
As is illustrated in Figure \ref{rad_curv_geod} down, $\vect n$ is a vector in $T_m\SS^n$ or $T_m\EE^n$ which is a subspace of $T_mH$. 
\begin{proposition}\label{prop_sigma=k_g m+n}
As above, the point $\sigma=\varphi(\Si)\in \Lambda$ is given by
\begin{equation}\label{f_sigma_kgmn}
\sigma = k_g m+\vect n.
\end{equation}
When the ambient space is Euclidean, the letter ``g" of $k_g$ can be dropped off.
\end{proposition}

\begin{remark}
If the orientation of $\Si$ is reversed then the corresponding point $\sigma$ in $\La$ should be replaced by $-\sigma$. 
Therefore, we need a sign convention for the geodesic curvature $k_g$, which we shall fix as follows. 
(We only explain in a spherical model.) 

Choose the unit normal vector $n$ so that if a basis of $T_m\Si$ consisting of ordered vectors $\vect v_1, \cdots, \vect v_{n-1}$ gives the positive orientation of $T_m\Si$, then a basis of $T_m\SS^n$ consisting of ordered vectors $\vect v_1, \cdots, \vect v_{n-1}, \vect n$ gives the positive orientation of $T_m\SS^n$. 
Let $\vect a$ be an accelaration vector at $m$ to a geodesic circle of $\Si$ through $m$ and $p$ the orthogonal projection to $T_m\SS^n$. 
Then $k_g$ is given by $p(\vect a)=k_g \vect n$. 
\end{remark}

\begin{proof} \rm 
We give two kinds of proofs in the case when $n=2$. 

\smallskip
(1) Assume $k_g\ne0$. 
First recall that if $q$ is a point in the light cone, the orthogonal complement of $\span{q}$ in the Minkowski space is the codimension $1$ hyperplane $T_q\lc$ which is tangent to the light cone along the line $\span{q}$. 
It implies that the line $\Pi^{\perp}$ which is orthogonal to $\Pi$ is contained in (in fact, it turns out to be equal to) $\displaystyle \bigcap_{q\in\Si}T_q\lc.$ 
By taking the intersection with the affine $3$-space $H$ (Figure \ref{rad_curv_geod}), it follows that $\Pi^{\perp}\cap H$ is the vertex $v$ of a cone which is tangent to $\SS^2(1)=\lc\cap H$ along $\Si$. 
Therefore $\sigma$ is given by $\sigma=\pm v/\Vert v\Vert$ (\cite{HJ}). 

Let $V$ be the cone with vertex $v$ which is tangent to $\SS^2$ (or $\EE^2$ or $\HH^2$) along $\Si$ (Figure \ref{rad_curv_geod} down). 
Develop the cone $V$ on an Euclidean plane by rolling it. 
As the cone is a developable ruled surface, this developing is a isometry. 
The curve obtained from $\Si$ is an arc of a circle whose curvature is equal to the geodesic curvature of $\Si$, which means that $\Vert v-m\Vert$ is equal to $1/k_g$ (Figure \ref{parabola_light_lines}). 
\begin{figure}[htbp]
\begin{center}
\includegraphics[width=.4\linewidth]{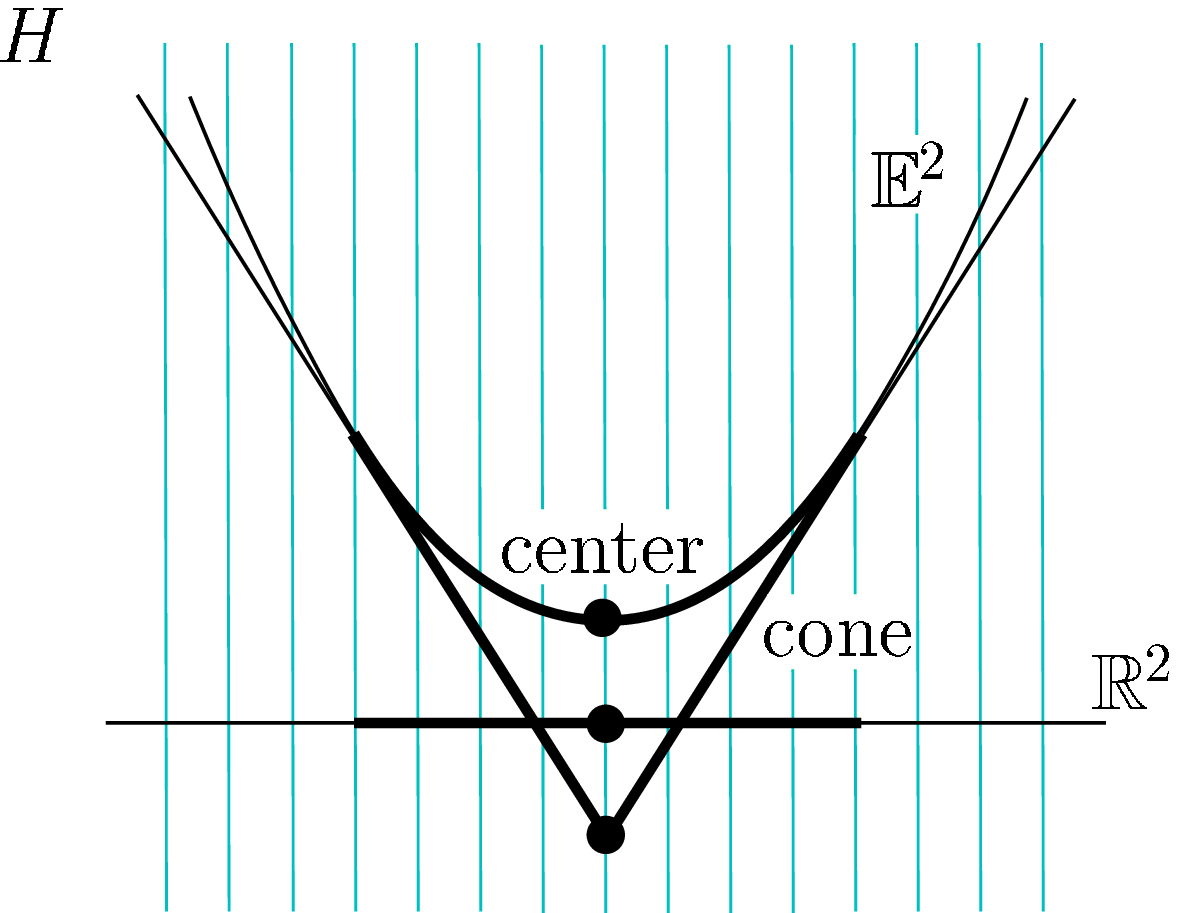}
\caption{The projection in the lightlike direction gives an isometry between $\EE^2$, $\RR^2$ and the cone in $H$ which is tangent to $\EE^2$ along $\Si$}
\label{parabola_light_lines}
\end{center}
\end{figure}

Therefore we have $\displaystyle v=m+\frac1{k_g}\vect n.$ 
Since $m$ is lightlike and orthogonal to $\vect n$, we have $\langle v,v\rangle=k_g{}^{-2}$, 
which implies that the unit normal vector to $\Pi$ is given by $\pm k_gv=k_gm+\vect n$. 

\smallskip
(2) We first give a proof for the spherical case. 
Let $\Si$ be a circle given by $\Si=\SS^2\cap W$ and $m$ a point in $\Si$. 
Since the bijection $\varphi$ is $O(4,1)$-equivariant (\ref{conf_equivar_varphi}), we may assume, after an action of $O(4,1)$ if necessary, that the spherical model is $\SS^2(1)$, and $W$ and $m$ are given by 
\setlength\arraycolsep{1pt}
\[\begin{array}{rcl}
W&=&\{(1,\cos\alpha,y,z)\,|\,y,z\in\mathbb R\}\hspace{0.5cm}(0<\alpha\le\frac{\pi}2), \\[1mm]
m&=& (1,\cos\alpha,\sin\alpha,0). \end{array}\]

Since the radius of $\Si$ in the affine space $W$ is equal to $\sin\alpha$ (Figure \ref{k_gm+n}) we have $\vect a=\displaystyle \left(0,0,-\frac1{\sin\alpha},0\right)$, 
\begin{figure}[htbp]
\begin{center}
\begin{minipage}{.45\linewidth}
\begin{center}
\includegraphics[width=.9\linewidth]{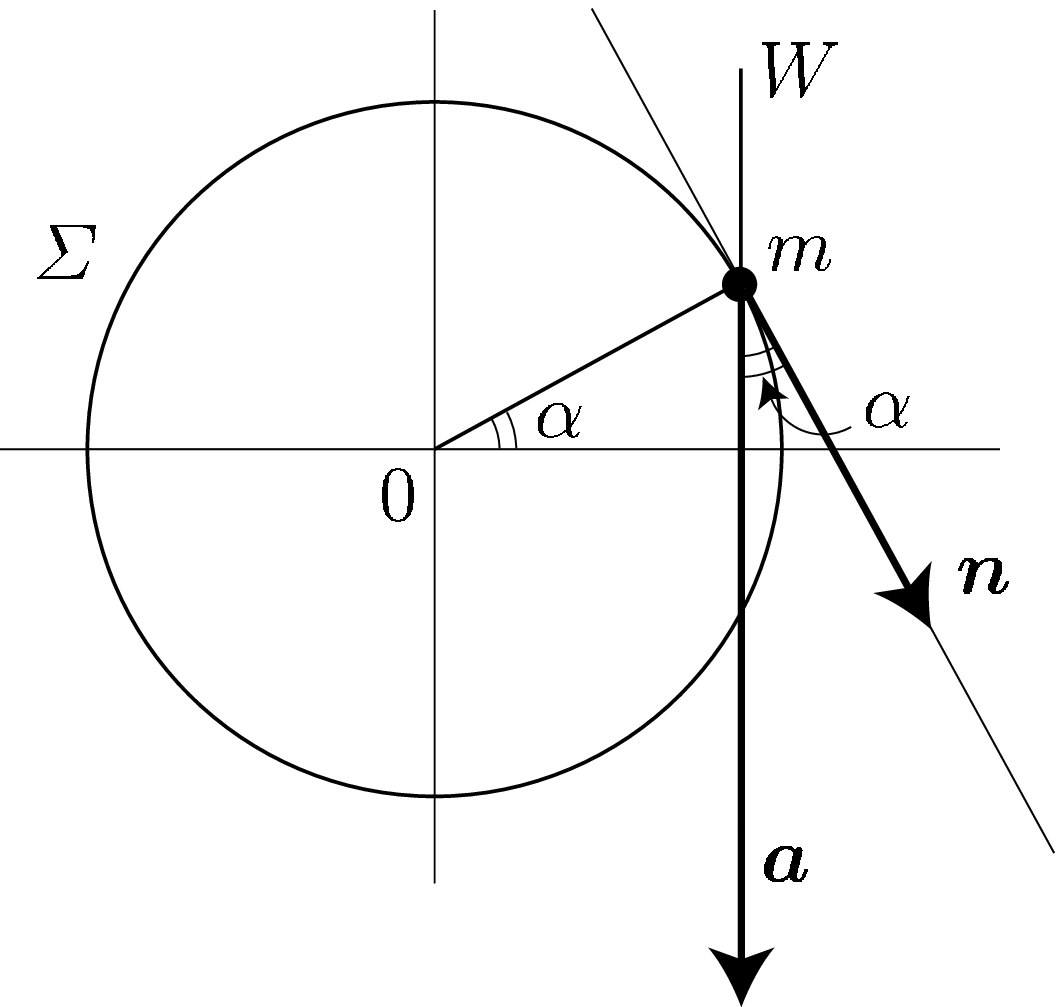}
\caption{A picture in the section with the plane containing the origin, $m$, and $\vect a$}
\label{k_gm+n}
\end{center}
\end{minipage}
\hskip 0.4cm
\begin{minipage}{.45\linewidth}
\begin{center}
\includegraphics[width=.9\linewidth]{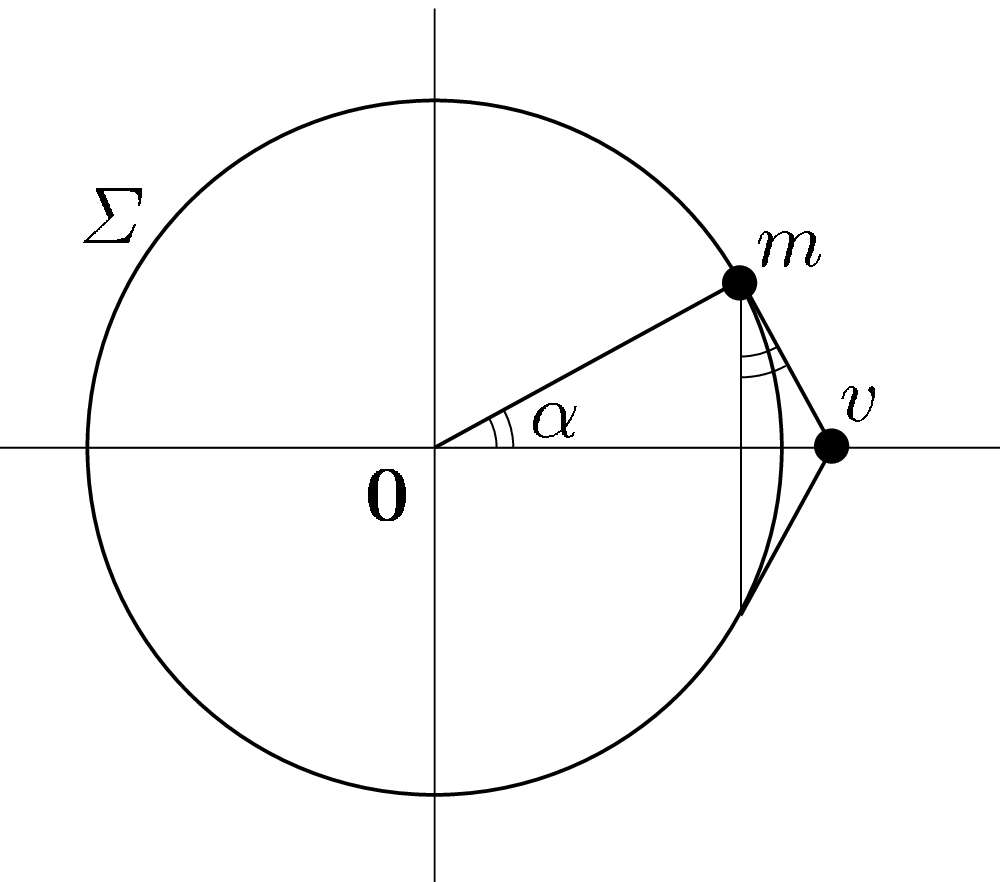}
\caption{The vertex $v$ of a cone tangent to $\SS^2(1)$ at $\Si$}
\label{m_v}
\end{center}
\end{minipage}
\end{center}
\end{figure}
which implies $p(\vect a)=\displaystyle \left(0,\cos\alpha,-\frac{\cos^2\alpha}{\sin\alpha},0\right)$. 
Since $\vect n=\pm(0,\sin\alpha,-\cos\alpha,0)$ it follows that $k_g=\pm\cot\alpha$. 

Since $\Pi$, where $\Si=\SS^2(1)\cap\Pi$, is given by $\Pi=\{c(1,\cos\alpha,y,z)\,|\,c,y,z\in\mathbb R\}$, its unit normal vector $\sigma$ is given by $\displaystyle \sigma=\pm\left(\cot\alpha,\frac1{\sin\alpha},0,0\right)$. 
Therefore we have $\sigma = k_g m+\vect n$. 

\smallskip
A proof for an the Euclidean case can be obtained exactly in the same way. 
The $O(4,1)$-equivariance of the bijection $\varphi$ enables us to locate the model $\EE^2$, the sphere $\Si$, and the point $m$ in special positions which make the computation very easy. 
\end{proof}
%

\section{{\boldmath$\frac12$} dimensional measure of lightlike curves}\label{sec_half_dim_length}
Let $\gamma$ be a lightlike curve in Minkowski space $\RR^{m}_1$, or in general, in a pseudo-Riemannian space $\RR^m_l$ with index $l$. 
(When $l>1$ a lightlike curve is called a {\em null} curve in \cite{O'N}. ) 
We define the Lorentz quadratic form $\LL(\vect v)$ by $\langle \vect v, \vect v\rangle$ and the norm by $\|v\|= \sqrt{|\LL(\vect v)|}\,$ for a vector $v$ in $\RR^m_l$ $(l>1)$ as well. 

\begin{proposition}
Let $\gamma$ be a compact lightlike curve which is piecewisely of class $C^4$. 
Let $t$ $(t\in[0,T])$ be the parameter, which is not necessarily the arc-length. 
We consider a subdivision $0=t_0 <t_1 \cdots <t_n = T$ of the interval $[0,T]$. 
Then the following limit exists, is finite, and generically non-zero: 
\begin{equation}\label{f_L1/2-measure}
{\lim_{\delta \rightarrow 0} \,\sum_i \sqrt{\| \gamma(t_{i+1})- \gamma (t_i) \|}\,,}
\end{equation}
where $\delta =  \max \{|t_{i+1}-t_i|\}$.
\end{proposition}
\begin{proof}
We have 
\[\gamma(t+h)-\gamma(t)=\Pt{\gamma}h+\frac{\Ppt{\gamma}}2h^2+\frac{\Pppt{\gamma}}6h^3+O(h^4).\]
Since $\gamma$ is lightlike, we have 
\[\big\langle \Pt{\gamma}, \>\Pt{\gamma}\big\rangle=0, \big\langle \Pt{\gamma}, \>\Ppt{\gamma}\big\rangle=0, \>\>\>\mbox{and}\>\>\>\big\langle \Ppt{\gamma}, \Ppt{\gamma}\big\rangle+\big\langle \Pt{\gamma}, \Pppt{\gamma}\big\rangle=0, \]
which implies

\setlength\arraycolsep{1pt}
\begin{equation}\label{1/12*ddot}
\begin{array}{rcl}
\LL(\gamma(t+h)-\gamma(t))&=&\langle \gamma(t+h)-\gamma(t), \gamma(t+h)-\gamma(t) \rangle \\[1mm]
&=&\displaystyle \left(\frac{\big\langle\Ppt{\gamma}, \Ppt{\gamma}\big\rangle}4+\frac{\big\langle\Pt{\gamma}, \Pppt{\gamma}\big\rangle}{3}\right)h^4+O(h^5)\\[2.5mm]
&=&\displaystyle -\frac1{12}\big\langle \Ppt{\gamma}, \Ppt{\gamma} \big\rangle \,h^4+O(h^5).
\end{array}
\end{equation}
It follows that 
\begin{equation}\label{f_L^(1/2)-measure_gamma_gamma''}
\lim_{\max|t_{j+1}-t_j|\to+0}\sum_{i}\sqrt{\Vert\gamma(t_{i+1})-\gamma(t_i)\Vert}=\int_C\sqrt[4]{\frac{|\LL(\Ppt{\gamma})|}{12}} \,dt\,. \end{equation}
\end{proof}
\begin{definition} \rm 
Let us call (\ref{f_L1/2-measure}) the {\em $\frac12$ dimensional measure} or {\em $L^{\frac12}$-measure} of a lightlike curve $\gamma$ and denote it by $L^{\frac12}(\gamma)$. 
\end{definition}
\begin{definition}\label{def_1/2-length_element} \rm 
Let $\gamma$ be a lightlike curve. 
Define a $1$-form $d\rho_{\!L^{\frac12}(\gamma)}$ on $\gamma$ by 
\begin{equation}\label{f_1/2-length_element}
\gamma^{\ast}d\rho_{\!L^{\frac12}(\gamma)}=\sqrt[4]{\frac{|\LL(\Ppt{\gamma})|}{12}} \,\,dt
\end{equation}
and call it the {\em $\frac12$ dimensional length element} or {\em $L^{\frac12}$-length element} of $\gamma$. 
\end{definition}
The formula (\ref{f_L^(1/2)-measure_gamma_gamma''}) implies that the $L^{\frac12}$-measure of a lightlike curve $\gamma$ satisfies
\begin{equation}\label{1/2-measure<->1/2-element}
L^{\frac12}(\gamma)=\int_{\gamma}d\rho_{\!L^{\frac12}(\gamma)}. 
\end{equation}

\begin{lemma}\label{lem_inv_form_null_curve} 
The $L^{\frac12}$-length element of a lightlike curve is well-defined, i.e. the right hand side of {\rm (\ref{f_1/2-length_element})} does not depend on the parametrization of $\gamma$ up to sign. 
\end{lemma} 
\begin{proof}
Let $t$ and $u$ be any parameters on $\gamma$. 
Let us denote $\frac{d}{du}$ by putting ${}^{\p}$ and $\frac{d}{dt}$ by $\Pt{}$. 
Then we have \setlength\arraycolsep{1pt}
\[\begin{array}{rcl}
\Pt{\gamma}&=&\displaystyle \frac{du}{dt}\,\gamma^{\p}\\[2mm]
\Ppt{\gamma}&=&\displaystyle \frac{d}{dt}\,\Pt{\gamma}=\frac{du}{dt}\cdot\frac{d}{du}\left(\frac{du}{dt}\,\gamma^{\p}\right)
= \frac{du}{dt}\left\{{\left(\frac{du}{dt}\right)}^{\p}\gamma^{\p}+\frac{du}{dt}\gamma^{\p\p}\right\}.
\end{array}
\]\setlength\arraycolsep{5pt}
Since $\langle \gamma^{\p}, \gamma^{\p}\rangle=0$ and hence $\langle \gamma^{\p}, \gamma^{\p\p}\rangle=0$ we have 
\[\left\langle\frac{d^2\gamma}{dt^2}, \frac{d^2\gamma}{dt^2}\right\rangle=\left(\frac{du}{dt}\right)^{\!4}\left\langle\frac{d^2\gamma}{du^2}, \frac{d^2\gamma}{du^2}\right\rangle,\]
which implies 
\[\sqrt[{4}]{\left|\left\langle\frac{d^2\gamma}{du^2}, \frac{d^2\gamma}{du^2}\right\rangle\right|}\,du
=\pm\,\sqrt[4]{\left|\left\langle\frac{d^2\gamma}{dt^2}, \frac{d^2\gamma}{dt^2}\right\rangle\right|}\,dt.\]
\end{proof}

\section{Lightlike curves in de Sitter spaces}\label{sec_lightlike_curves}

In this section we consider two examples of lightlike curves in de Sitter space: a curve in $\Lambda^3$ consisting of the osculating circles to a curve in $\EE^2$ or $\SS^2$, and that in $\Lambda^4$ of the focal spheres to a surface in $\EE^3$ or $\SS^3$ along a corresponding line of principal curvature of the surface. 

\subsection{Lightlike curves in $\mathcal{S}(1,2)$}\label{subs_light_12}

The set $\mathcal{S}(1,2)$ of oriented circles in $\SS^2$ can be identified with de Sitter space $\La^3$ in $\RR^4_1$. 
We give a characterization of a set of osculating circles to a curve in $\SS^2$ or $\EE^2$. 

Recall that a point in a plane curve is a {\em vertex} if and only if $k^{\p}=0$ holds at that point, which happens if and only if the curve has the third order contact at that point. 
As the second condition is conformally invariant, we adopt it as the definition of a vertex of a curve in a sphere (Definition \ref{def_vertex}). 

\begin{theorem}\label{lem_lightlike_curve=osc.cirlces}
A curve $\gamma$ in $\mathcal{S}(1,2)$ is a set of the osculating circles to a {vertex-free} curve $C$ in $\SS^2$ or $\EE^2$ with a non-vanishing velocity vector if and only if $\gamma$ is lightlike {and $\dim\span{\Pt{\gamma}, \Ppt{\gamma}}=2$}. 

Without the two conditions, $C$ being vertex-free and $\dim\span{\Pt{\gamma}, \Ppt{\gamma}}=2$, the ``if'' part of the above statement may fail although the ``only if'' part still holds. 
\end{theorem}

\begin{proof} 
We prove the first statement first. 

(1-1) ``Only if\,'' part. 
Suppose $\gamma$ is the set of osculating circles to a curve $C=\{m(s)\}$, where $s$ is the arc-length of $C$. 

Let $T$ denote the unit tangent vector to $C$: $T={m}^{\p}$. 
The osculating circle to the curve $C$ at a point $m(s)$ has the same geodesic curvature as that of $C$ at the same point, which is $k_g(s)$.  
Therefore Proposition \ref{prop_sigma=k_g m+n} implies that $\gamma(s)$ is given by 
$\gamma(s)= k_g(s) m(s)+\vect n(s),$ 
which implies that 
\[\gamma^{\p}(s)={k_g}^{\p}(s)m(s) + k_g(s) {m}^{\p}(s) + {\vect n}^{\p}(s).\] 
As ${m}^{\p}=T$ and ${\vect n}^{\p}= -k_g T$, we have ${\gamma}^{\p}(s)= {k_g}^{\p}(s)m(s)$, which proves that $\gamma$ is lightlike.

Note that $m={\gamma}^{\p}/{k_g}$ because ${k_g}^{\p}\ne 0$ as $C$ is vertex-free. 
Since $m(s)$ is in an affine subspace $H$ that does not passes through the origin, $\dim\langle{m, {m}^{\p}}\rangle=2$, which implies $\dim\langle{{\gamma}^{\p}, {\gamma}^{\p\p}}\rangle=2$. 

\smallskip
(1-2) ``If\,'' part. 
Suppose $\gamma$ is a lightlike curve with non-vanishing tangent vectors. 
Each tangent vector defines a point in $\SS^2$ by $\SS^2\cap\span{\Pt{\gamma}(t)}$, which we denote by $f(t)$. 
Put $C=\{f(t)\}$. 

Let us denote the derivation with respect to $t$ by putting $\Pt{}$ above. 
We may assume, after reparametrization if necessary, that the $x_0$-coordinate of $\Pt{\gamma}(t)$ is always equal to $1$. 
Then the point $f(t)$ in $\SS^2(1)$ is given by $f(t)=\Pt{\gamma}(t)$. 

As $\langle\gamma(t), \gamma(t)\rangle=1$ and $\langle\Pt{\gamma}(t), \Pt{\gamma}(t)\rangle=0$ we have 
\[\langle\gamma(t), \Pt{\gamma}(t)\rangle=\langle\gamma(t), \Ppt{\gamma}(t)\rangle=\langle\gamma(t), \Pppt{\gamma}(t)\rangle=0,\]
which implies 
\[\gamma\in\left(\spanbig{\Pt{\gamma}(t), \Ppt{\gamma}(t), \Pppt{\gamma}(t)}\right)^{\perp}
=\left(\vspace{-0.1cm}\spanbig{f(t), \Pt{f}(t), \Ppt{f}(t)}\right)^{\perp}. \]
Since $\dim\spanbig{f, \Pt{f}}=\dim\spanbig{\Pt{\gamma}, \Ppt{\gamma}}=2$ by the assumption, $\Pt{f}$ never vanishes, therefore 
$\dim \spanbig{f(t), \Pt{f}(t), \Ppt{f}(t)}=3.$ 
Since the osculating circle to $C$ at $f(t)$ is given by $\SS^2\cap\spanbig{f(t), \Pt{f}(t), \Ppt{f}(t)}$, it corresponds to $\pm \gamma(t)$ in $\La^3$. 

\smallskip
(2) Suppose $\gamma$ is a lightlike line in de Sitter space. 
Then it corresponds to a family of circles which are all tangent to each other at a constant point, which cannot be a family of osculating circles to a curve in $\SS^2$ or $\EE^2$. 
As $\Ppt{\gamma}\equiv \vect 0$ in this case, it implies the second statement of the theorem. 
\end{proof}
Let us prove that the $L^{\frac{1}{2}}$-measure of the lightlike curve $\gamma \subset \Lambda^3$ is equal to the conformal arc-length of $C$.
\begin{figure}[htb]
\begin{center}
\includegraphics[width=35mm]{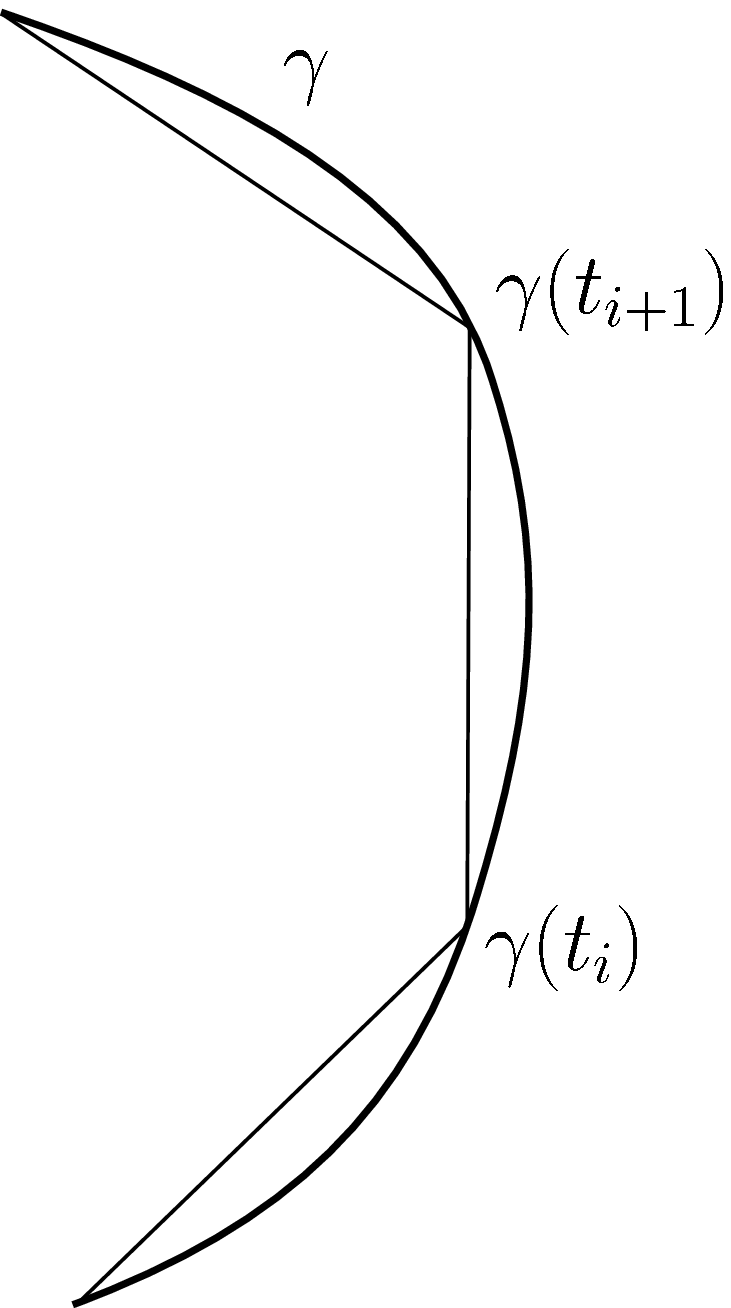}
\caption{An inscribed polygon on a lightlike curve $\gamma$. A vector $\gamma(t_{i+1})-\gamma(t_i)$ is timelike in general.}
\label{length_light}
\end{center}
\end{figure}

Let $s$ be the arc-length of a curve $C$. 
\begin{theorem}
The conformal length $\rho$ of a curve $C$ in $\ss^2$, $\RR^2$, or $\HH^2$ is equal to $\sqrt[4]{12}$ times the $\frac{1}{2}$-dimensional measure $L^{\frac{1}{2}}(\gamma)$ of the lightlike curve $\gamma \subset \Lambda^3$ which consists of the osculating circles to $C$. 
It is given by $\displaystyle \int_C \sqrt{|k'_g|} ds$ (we can drop the letter ``g" when $C$ is a plane curve). 
\end{theorem}
\begin{proof}\label{length_osc_circles}
Let us denote the derivation with respect to $s$ by putting ${\phantom{|}}^{\p}$ as before. 
Proposition \ref{prop_sigma=k_g m+n} and the proof of ``only if\,'' part of Theorem \ref{lem_lightlike_curve=osc.cirlces} imply that $\gamma=k_g m+n$ and $\gamma^{\p}={k_g}^{\p}m$. 
Therefore the second derivative ${\gamma}^{\p\p}$ is given by ${k_g}^{\p\p}m + {k_g}^{\p}T$, where $T=m^{\p}$ is the unit tangent vector to $C$. 
As the unit vector $T$ is orthogonal to the lightlike vector $m$ and spacelike since all the tangent vectors to the light cone which are not tangent to a generatrix are spacelike, we see that  
\[\LL(\gamma^{\p\p})=\LL({k_g}^{\p\p}m + {k_g}^{\p}T)={{k_g}^{\p}}^2\LL(T)= {{k_g}^{\p}}^2.\]
Therefore, (\ref{f_L^(1/2)-measure_gamma_gamma''}) implies that the $L^{\frac{1}{2}}$-measure of $\gamma$ is given by 
$$L^{\frac{1}{2}}(\gamma) = \sqrt[4]{\frac{1}{12}}\int_C \sqrt{|{k_g}^{\p}|}\,ds.$$
Since $\sqrt{|{k_g}^{\p}|}\,ds=d\rho$ by (\ref{inf_conf_arc-length}), it completes the proof. 
\end{proof}
Note that (\ref{1/12*ddot}) implies that when ${k_g}^{\p}\not= 0$ the vector $\gamma(t+h)-\gamma(t)$ is timelike  for small enough $h$. 
\subsection{Lightlike curves in $\mathcal{S}(2,3)$}\label{subs_light_23}
We say that a sphere $\Si_O$ is an {\em osculating sphere} of a curve $C=\{f(t)\}$ if it has the third order contact with $C$. 
It happens if and only if the point $\varphi(\Si_O)$ in $\La$ which corresponds to $\Si_O$ satisfies $\varphi(\Si_O)\perp\spanbig{f, \Pt{f}, \Ppt{f}, \Pppt{f}}$. 

\smallskip
\begin{remark}\rm 
We use the word ``osculating spheres'' for curves and ``focal spheres'' for surfaces (Theorem \ref{thm_focal_surfaces}). 
\end{remark}

\smallskip
We identify the space $\mathcal{S}(2,3)$ of oriented spheres in $\SS^3$ with the $4$-dimensional de Sitter space $\La$ in $\RR^5_1$ as before. 
\begin{proposition} 
Suppose $\sigma$ is a lightlike curve in $\mathcal{S}(2,3)$. 
Since $\sigma^{\p}(t)$ is lightlike, it defines a point in $\SS^3$ by $\SS^3\cap\span{\sigma^{\p}(t)}$, which shall be denoted by $f(t)$. 
Put $C=\{f(t)\}$. 
Then a sphere $\Si(t)$ which corresponds to $\sigma(t)$ is not necessarily an osculating sphere of $C$, but it contains the osculating circle to $C$ at $f(t)$. 
\end{proposition}
\begin{proof} 
The proof is parallel to that of the ``if\,'' part of Theorem \ref{lem_lightlike_curve=osc.cirlces}. 
As $\langle\sigma(t), \sigma(t)\rangle=1$ and $\langle\Pt{\sigma}(t), \Pt{\sigma}(t)\rangle=0$ we have 
\[\langle\sigma(t), \Pt{\sigma}(t)\rangle=\langle\sigma(t), \Ppt{\sigma}(t)\rangle=\langle\sigma(t), \Pppt{\sigma}(t)\rangle=0.\]
Therefore $\span{\sigma(t)}$ is contained in $\left(\span{\Pt{\sigma}(t), \Ppt{\sigma}(t), \Pppt{\sigma}(t)}\right)^{\perp}$, which implies $\Si(t)=\SS^3\cap(\span{\sigma(t)})^{\perp}$ contains $\SS^3\cap\span{\Pt{\sigma}(t), \Ppt{\sigma}(t), \Pppt{\sigma}(t)}$ which is the osculating circle of $C$ at $f(t)$. 
\end{proof}

\smallskip
Let us now consider a surface $M$ in $\ss^3$, $\HH^3$, or $\RR^3$.
In the Poincar\'e ball model, the intersection of a $2$-sphere with the ball is called a sphere. It may be a usual geodesic sphere or a plane of constant curvature between $0$ and $-1$.

At a point $m\in M$, a sphere tangent to $M$ at $m$ whose geodesic curvature is equal to one of the principal curvatures $k_1,k_2$ of $M$ at $m$, has higher contact with $M$. 
Let us denote them by $\Si_1$ and $\Si_2$ and call them {\em focal spheres} to $M$ at $m$. 
They are distinct if the point $m$ is not an umbilical point of $M$. 
\begin{theorem}\label{thm_focal_surfaces}
The curve $\gamma_1 \subset \Lambda^4$ corresponding to the focal spheres $\Si_1$ along a line of principal curvature $C_1$ for the principal curvature $k_1$ is lightlike. Its $\frac{1}{2}$-dimensional measure is 
$$L^{\frac{1}{2}}(\gamma_1)=\sqrt[4]{\frac{1}{12}} \int_{C_1} \sqrt{|X_1(k_1)|}\, ds$$
\noindent where $X_1$ is the unit tangent vector to $C_1$ which is parametrized by its arc-length $s$. 
\end{theorem}
\begin{proof}
The proof is the same as above (Theorem \ref{length_osc_circles}) using the formula $\gamma(s)=k_gm(s)+\vect n(s)$, where $k_g=k_1$ is the geodesic curvature of the focal sphere $\Si_1$. The hypothesis that $C_1$ is a line of principal curvature associated to the principal curvature $k_1$ implies that $\vect n^{\p}(s)=k_gX_1(s)$. 
\end{proof}
\begin{remark}
A similar statement is also valid if $M$ is an hypersurface of some space-form. The integral $\int_{C_1} \sqrt{|X_1(k_1)|} ds$ has already appeared in \cite{Ro-Sa}.
\end{remark}
\section{Space $\mathcal{S}(1,3)$ of the oriented circles in $\SS^3$}\label{sec_S(1,3)}
The pseudo-Riemannian structure of the space of the oriented circles in $\RR^3$ (or $\SS^3$) also arises naturally in the study of conformal geometry. 
Each tangent space of $\mathcal{S}(1,3)$ has an indefinite non-degenerate quadratic form which is compatible with M\"obius transformations.

\subsection{The set of circles as a Grassmann manifold}
An oriented circle in $\RR^3$ (or $\SS^3$) can be realized in the Minkowski space $\mathbb{R}^{5}_1$ as the intersection of the light cone and an oriented timelike $3$-dimensional vector subspace. 
Therefore the set $\mathcal{S}(1,3)$ can be identified with the Grassmann manifold $\widetilde{\textsl{Gr}}_-(3;\mathbb{R}^{5}_1)$ of oriented $3$-dimensional timelike subspaces of $\mathbb{R}^{5}_1$. 

It follows that the set $\mathcal{S}(1,3)$ is the homogeneous space $SO(4,1)/SO(2)\times SO(2,1).$ 
We give its pseudo-Riemannian structure explicitly in what follows. 

\subsection{Pl\"ucker coordinates for the set of circles}
Let us recall the Pl\"ucker coordinates of Grassmannian manifolds. 

Let $W$ be an oriented $3$-dimensional vector subspace in $\mathbb{R}^{5}_1$, and let $\{\vect x_1, \vect x_2, \vect x_{3}\}$ be an ordered basis of $W$ which gives the orientation of $W$. 
Define $p_{i_1i_2i_3}$ $(0\le i_k\le 4)$ by 
\begin{eqnarray}\label{f_def_plucker_coordinates}
p_{i_1i_2i_{3}}=
\left|\begin{array}{ccc}
x_{1\,i_1}& x_{1\,i_2} &x_{1\,i_{3}}\\
x_{2\,i_1}&x_{2\,i_2}&x_{2\,i_{3}}\\
x_{3\,i_1}&x_{3\,i_2}&x_{3\,i_{3}}
\end{array}\right|.
\end{eqnarray}

Let $[W]$ denote an unoriented $3$-space which is obtained from $W$ by forgetting its orientation. 
Then it can be identified by the homogeneous coordinates $[\cdots,p_{i_1i_2i_3},\cdots]\in\mathbb{R}P^{9}$ called the {\em Pl\"ucker coordinates} or {\em Grassmann coordinates}. 
They do not depend on the choice of a base of $[W]$. 

The Pl\"ucker coordinates $p_{i_1i_2i_3}$ are not independent. 
They satisfy the {\em Pl\"ucker relations}: 
\begin{equation}\label{plucker_relations}
\sum_{k=1}^{4}(-1)^{k}p_{i_1i_{2}j_k}p_{j_1\cdots \widehat{j_k}\cdots j_{4}}=0,
\end{equation}
where $\widehat{j_k}$ indicates that the index $j_k$ is being removed. 
There are five non-trivial Pl\"ucker relations and exactly three of them are independent. 

As we are concerned with the orientation of the subspaces, we use the Euclidean spaces for the Pl\"ucker coordinates in this article instead of the projective spaces which are used in most cases. 
The {\em exterior product} of $\vect x_1, \vect x_2$, and $\vect x_{3}$ in $\RR^5_1$ is given by
$$\vect x_1\w\vect x_2\w\vect x_{3}=(\cdots,p_{i_1i_2i_{3}},\cdots)\in\mathbb R^{10} \hspace{0.5cm} (i_1<i_2<i_3)$$
through the identification ${\stackrel{3}{\mbox{$\bigwedge$}}}\,\mathbb{R}^{5}\cong\RR^{10}$. 

Let $\widetilde{\textsl{Gr}}_-(3;\RR^5_1)$ denote the Grassmann manifold of the set of all oriented $3$-dimensional timelike vector subspaces in $\RR^5_1$. 
\subsection{Pseudo-Riemannian structure of $\displaystyle {\stackrel{3}{\mbox{$\bigwedge$}}}\,\mathbb{R}^{5}_1$} 
The indefinite inner product of the Minkowski space $\mathbb{R}^{5}_1$ naturally induces that of $\displaystyle {\stackrel{3}{\mbox{$\bigwedge$}}}\,\mathbb{R}^{5}_1$ by  
\begin{equation}\label{f_wedge-det_x}
\langle \vect x_1\w\vect x_2\w\vect x_3, \vect y_1\w\vect y_2\w\vect y_3\rangle
=-\left|\!
\begin{array}{ccc}
\langle \vect x_1, \vect y_1\rangle & \langle \vect x_1, \vect y_2\rangle & \langle \vect x_1, \vect y_3\rangle \\
\langle \vect x_2, \vect y_1\rangle & \langle \vect x_2, \vect y_2\rangle & \langle \vect x_2, \vect y_3\rangle \\
\langle \vect x_3, \vect y_1\rangle & \langle \vect x_3, \vect y_2\rangle & \langle \vect x_3, \vect y_3\rangle 
\end{array}\!
\right|
\end{equation}
for any $\vect x_i$ and $\vect y_j$ in $\RR^5_1$. 
Note that this sign convention is opposite to that in \cite{HJ}. 

The above formula is a natural generalization of the one for the exterior products of four vectors which we obtained in \cite{La-OH1}, where we studied the set of oriented spheres in $\SS^3$. 

It follows that $\displaystyle {\stackrel{3}{\mbox{$\bigwedge$}}}\,\mathbb{R}^{5}_1$ can be identified with $\mathbb R^{10}$ with a pseudo-Riemannian structure with index $4$, which we denote by $\mathbb R^{10}_4$, so that $\{\vect e_{i_1}\w\vect e_{i_2}\w\vect e_{i_{3}}\}_{i_1<i_2<i_{3}}$ is a pseudo-orthonormal basis of $\displaystyle {\stackrel{3}{\mbox{$\bigwedge$}}}\,\mathbb{R}^{5}_1$ with 
\begin{equation}\label{f_orth_basis_ext_p}
\langle \vect e_{i_1}\w\vect e_{i_2}\w\vect e_{i_{3}}, \,\vect e_{i_1}\w\vect e_{i_2}\w\vect e_{i_{3}}\rangle=\left\{
\begin{array}{lcl}

-1 & \textrm{ if } & i_1\ge 1,\\[1mm]
+1 & \textrm{ if } & i_1=0.
\end{array}
\right.
\end{equation}
Let us realize $\mathcal{S}(1,3)$ as a pseudo-Riemannian submanifold of $\displaystyle {\stackrel{3}{\mbox{$\bigwedge$}}}\,\mathbb{R}^{5}_1\cong\mathbb R^{10}_{4}$. 

\begin{lemma}\label{lemm_time_cond}
Let $W$ be an oriented $3$-dimensional vector subspace in $\mathbb{R}^{5}_1$ spanned by $\vect x_1, \vect x_2, \vect x_{3}$. 
Then $W$ is timelike if and only if 
\[\langle \vect x_1\w\vect x_2\w\vect x_{3}, \,\vect x_1\w\vect x_2\w\vect x_{3}\rangle>0,\] 
and isotropic (i.e. tangent to the light cone) if and only if 
\[\langle \vect x_1\w\vect x_2\w\vect x_{3}, \,\vect x_1\w\vect x_2\w\vect x_{3}\rangle=0.\] 
\end{lemma}

\begin{proof}
Case (1). Suppose $W$ is not isotropic. 
We may assume without loss of generality that $\{\vect x_1, \vect x_2, \vect x_3\}$ is a pseudo-orthonormal basis of $W$. 
If $W$ is timelike, one of $\vect x_1, \vect x_2$ and $\vect x_3$ is timelike, and therefore $\langle \vect x_1\w\vect x_2\w\vect x_{3}, \,\vect x_1\w\vect x_2\w\vect x_{3}\rangle=1$ by (\ref{f_wedge-det_x}). 
If $W$ is spacelike, then $\langle \vect x_1\w\vect x_2\w\vect x_{3}, \,\vect x_1\w\vect x_2\w\vect x_{3}\rangle=-1$. 

Case (2). Suppose $W$ is isotropic. 
Then $W$ is tangent to the light cone at a lightlike line $l$. 
Now we may assume without loss of generality that $\{\vect x_1, \vect x_2, \vect x_3\}$ is a pseudo-orthogonal basis of $W$ and that $\vect x_1$ belongs to $l$. 
Then we have  $\langle \vect x_1\w\vect x_2\w\vect x_{3}, \,\vect x_1\w\vect x_2\w\vect x_{3}\rangle=0$. 
\end{proof}

\subsection{Conformal invariance of the pseudo-Riemannian structure} 
We show that $\mathcal{S}(1,3)$ has a pseudo-Riemannian structure which is compatible with M\"obius transformations of $\SS^3$. 
Let $O(6,4)$ denote the pseudoorthogonal group of ${\stackrel{3}{\mbox{$\bigwedge$}}} \,\mathbb{R}^{5}_1\cong\RR^{10}_4$. 
\begin{definition}\label{def_Psi_A} \rm 
Define a map $\varPsi:M_{5}(\mathbb{R})\to M_{10}(\mathbb{R})$ by 
$$
\varPsi:M_{5}(\mathbb{R})\ni A=(a_{ij})\mapsto \varPsi(A)=(\tilde a_{IJ})\in M_{10}(\mathbb{R}),
$$
where $I=(i_1i_2i_{3})$ and $J=(j_1j_2j_{3})$ are multi-indices, 
and $\tilde a_{IJ}$ is given by 
$$
\tilde a_{IJ}=
\left|\!
\begin{array}{ccc}
a_{i_{1}j_{1}}& a_{i_{1}j_{2}} & a_{i_{1}j_{3}}\\
a_{i_{2}j_{1}}& a_{i_{2}j_{2}} & a_{i_{2}j_{3}} \\
a_{i_{3}j_{1}}& a_{i_{3}j_{2}} & a_{i_{3}j_{3}}
\end{array}\!
\right|. 
$$
\end{definition}

\begin{lemma}\label{prop_tilde_A} 
\begin{enumerate}
\item We have 
\begin{equation}\label{charact_Psi_A}
(A\vect x_1)\w (A\vect x_2) \w(A\vect x_{3})=\varPsi(A)\,(\vect x_1\w \vect x_2 \w \vect x_{3}) \hspace{0.5cm} (\forall \vect x_1, \vect x_2, \vect x_{3}\in\mathbb{R}^{5}_1) 
\end{equation}
for $A\in M_5(\RR)$. 
\item If $A\in O(4,1)$ then $\varPsi(A)\in O(6,4)$. 
\item The restriction of $\varPsi$ to $O(4,1)$ is a homomorphism. 
\end{enumerate}
\end{lemma}

We can say more, although we do not give proof: 
The matrix $\varPsi(A)$ can be characterized by (\ref{charact_Psi_A}). 
The reverse statement of (2) also holds. 
The restriction of $\varPsi$ to $\textsl{Gl}(5,\RR)$ is a homomorphism whose kernel consists of $\{\pm I\}$. 
\begin{proof}
(1) The definition of $\tilde a_{IJ}$ implies 
\[\varPsi(A)\,(\vect e_{i_1}\w \vect e_{i_2}\w \vect e_{i_3})=(A\vect e_{i_1})\w(A\vect e_{i_2})\w(A\vect e_{i_3}).\]

(2) If $A\in O(4,1)$ then (\ref{f_wedge-det_x}) and (\ref{charact_Psi_A}) imply 
\[\left\langle \varPsi(A)(\vect e_{i_1}\w\vect e_{i_2}\w\vect e_{i_{3}}), \,\varPsi(A)(\vect e_{j_1}\w\vect e_{j_2}\w\vect e_{j_{3}}) \right\rangle=\left\langle \vect e_{i_1}\w\vect e_{i_{2}}\w\vect e_{i_{3}}, \,\vect e_{j_1}\w\vect e_{j_{2}}\w\vect e_{j_{3}} \right\rangle,\]
which implies $\varPsi(A)\in O(6,4)$. 

(3) Routine calculation in linear algebra implies $\varPsi(AB)_{IJ}=\sum_K\tilde a_{IK}\tilde b_{KJ}$. 
\end{proof}

\begin{corollary}
We have 
$$\psi(A\cdot\Ga)=\varPsi(A)\psi(\Ga)$$
for $\Ga\in\mathcal{S}(1,3)$ and  $A\in O(4,1)$, where $\psi$ is the bijection from $\mathcal{S}(1,3)$ to $\Theta(1,3)\subset\mathbb R^{10}_{4}$ given by {\rm (\ref{S=Gr=Theta})} and $\varPsi$ the homomorphism from $O(4,1)$ to $O(6,4)$ given in {\rm Definition \ref{def_Psi_A}}. 
\end{corollary}

\begin{proposition}\label{cor_S=Theta-pR-str}  
Let $\Theta(1,3)$ be the intersection of the quadric satisfying the Pl\"ucker relations and the unit pseudo-sphere: 
\begin{equation}\label{def_Theta}
\Theta(1,3) = \left\{
(\cdots,p_{i_1i_2i_{3}},\cdots)
\in\mathbb R^{10}_{4}
\left|\begin{array}{l}
\displaystyle \sum_{k=1}^{4}(-1)^{k}p_{i_1i_{2}j_k}p_{j_1\cdots \widehat{j_k}\cdots j_{4}}=0\\[3.5mm]
\displaystyle -\sum_{i_1\ge1}p_{i_1i_2i_{3}}{}^2+\sum_{i_2\ge 1}p_{0i_2i_{3}}{}^2=1
\end{array}
\right\}\right..
\end{equation}
It is a $6$-dimensional pseudo-Riemannian submanifold of $\mathbb R^{10}_{4}$ with pseudo-Riemannian structure of index $2$. 

Then the set $\mathcal{S}(1,3)$ of oriented circles in $\SS^3$ can be identified with $\Theta(1,3)$ through a bijection $\psi$ given by 
\begin{equation}\label{S=Gr=Theta}
\begin{array}{ccccccc}
\psi:\hspace{-0.4cm}&\!\!\mathcal{S}(1,3)\!\!&\!\!\stackrel{\displaystyle \cong}{\longrightarrow}\!\!&\!\!\widetilde{\textsl{Gr}}_-(3;\mathbb{R}^{5}_1)\!\!&\!\!\stackrel{\displaystyle \cong}{\longrightarrow}\!\!&\!\!\Theta(1,3)&\hspace{-0.8cm}\subset {\stackrel{3}{\mbox{$\bigwedge$}}}\,\mathbb{R}^{5}_1\cong \mathbb R^{10}_{4}\\[1mm] 
\!\!&\!\!\mbox{\rotatebox{90}{$\in$}}&&\mbox{\rotatebox{90}{$\in$}}&&\mbox{\rotatebox{90}{$\in$}}&\\[1mm]
\!\!&\!\!W\cap \SS^3(\infty) \!\!&\!\! \mapsto \!\!&\!\! W=\textsl{Span}\langle\vect x_1, \vect x_2, \vect x_{3}\rangle \!\!&\!\! \mapsto \!\!&\!\! \displaystyle \frac{\vect x_1\w\vect x_2\w\vect x_{3}}{\,\|\vect x_1\w\vect x_2\w\vect x_{3}\|\,}\,.&
\end{array}
\end{equation}
\end{proposition}

\begin{proof}
Since $\mathcal{S}(1,3)$ can be identified with the Grassmann manifold $\widetilde{\textsl{Gr}}_-(3;\mathbb{R}^{5}_1)$ of oriented $3$-dimensional timelike subspaces of $\mathbb{R}^{5}_1$, Lemma \ref{lemm_time_cond} implies that it is enough to show that the restriction of the indefinite inner product of $\mathbb R^{10}_{4}$ to each tangent space of $\Theta(1,3)$ induces a non-degenerate quadratic form of index $2$. 

The conformal invariance of the pseudo-Riemannian structure allows us to assume that an oriented circle $\Ga$ passes through $(\pm 1,0,0,0)$ and $(0,1,0,0)$. 
The index can be calculated in several ways. 

(i) $\Ga$ corresponds to $W=\textsl{Span}\langle\vect e_0, \vect e_1, \vect e_2\rangle$ in the Grassmannian $\widetilde{\textsl{Gr}}_-(3;\mathbb{R}^{5}_1)$. 
The tangent space $T_{W}\widetilde{\textsl{Gr}}_-(2;\mathbb{R}^{5}_1)$ is isomorphic to $\textrm{Hom}\left(W, W^{\perp}\right)$, which is isomorphic to $M_{3,2}(\RR)$. 
We can construct six vectors which form a pseudo-\linebreak orthonormal basis of the tangent space explicitly. 
It turns out that two of them are timelike and the other four are spacelike. 

(ii) The tangent space $T_{\Ga}\Theta(1,3)$ can be identified with the pseudo-orthogonal complements of the subspace spanned by gradients of the defining functions of $\Theta(1,3)$ which appear in (\ref{def_Theta}). 
There are five non-trivial Pl\"ucker relations and exactly three of them are independent.
Two of them give timelike gradients and the rest gives a spacelike one. 
On the other hand, the gradient of $-\sum_{i_1\ge1}p_{i_1i_2}{}^2+\sum_{i_2\ge 1}p_{0i_2}{}^2-1$ is spacelike. 
Hence the index can be given by $4-2=2$. 
\end{proof}

\medskip
Since $\mathcal{S}(1,3)$ is the homogeneous space $SO(4,1)/SO(2)\times SO(2,1)$, Proposition 3.2.6 of \cite{Ko-Yo} also implies that the index of $\Theta(1,3)$ is equal to $2$.  

\section{Osculating circles and the conformal arc-length}\label{sec_osculating_circles}
Let us realize the Euclidean space $\RR^3$ in the Minkowski space $\RR^5_1$ as the isotropic affine section of the light cone $\lc$ given by (\ref{def_E^3_0}) in section 
\ref{sec_preliminaries}. 
We use the following notation in what follows. 
Let $C=\{m(s)\}$ be an oriented curve in $\RR^3$ parametrized by the arc-length $s$. 
Let $\bar{m}$ be a map which is induced from $m$; 
\begin{equation}\label{m_bar_m}
\bar{m}(s)=\left(1+\frac{m(s)\cdot m(s)}4,\, -1+\frac{m(s)\cdot m(s)}4,\, m(s)\right)\in\EE^3_0\subset\RR^5_1,
\end{equation}
where $\EE^3_0$ is given by formula (\ref{def_E^3_0}). 

\smallskip
The {\em osculating circle} of a curve $C$ at a point $x$ is the circle with the best contact with $C$ at $x$. 
We will denote it by ${\mathcal O}_{x}$. 
It has the second order contact with $C$ at $x$. 

Suppose $x$, $y$, and $z$ are points on $C$. 
When $x$, $y$, and $z$ are mutually distinct, let $\Ga(x,y,z)$ denote the circle that passes through the points $x$, $y$, and $z$ in $C$ whose orientation is given by the cyclic order of $\{x,y,z\}$. 
When two (or three) of the points $x$, $y$, and $z$ coincide, $\Ga(x,y,z)$ means a tangent circle (or respectively, an osculating circle) whose orientation coincides with that of $C$ at the tangent point. 

\smallskip
Let $\gamma(u,v,w)$ be a point in $\Theta(1,3)\subset \mathbb R^{10}_4$ which corresponds to a circle $\Ga(\bar{m}(u), \bar{m}(v), \bar{m}(w))$ 
through the bijection $\psi$ from $\mathcal{S}(1,3)$ to $\Theta(1,3)$ (see (\ref{S=Gr=Theta})). 
Put $\gamma(u)=\gamma(u,u,u)$. 
It corresponds to the osculating circle at ${m}(u)$. 

\subsection{The curve of the osculating circles is lightlike}
Let us express the osculating circles using the exterior products of vectors in the Minkowski space. 

Observe that $\langle \bar{m}, \bar{m}\rangle \equiv 0$, 
as $\bar{m}$ belongs to the light cone, and that 
\[\langle {\bar{m}\phantom{|}\!}^{\p}, {\bar{m}\phantom{|}\!}^{\p}\rangle \equiv  {{m}}^{\p}\cdot{{m}}^{\p}\equiv 1,\] 
where $\cdot$ denotes the standard inner product of $\mathbb R^3$. 
Define $F_2$ and $F_3$  by 
\begin{equation}\label{def_F2F3}
F_2=\langle {\bar{m}\phantom{|}\!}^{\p\p}, {\bar{m}\phantom{|}\!}^{\p\p}\rangle, \>
F_3=\langle {\bar{m}\phantom{|}\!}^{\p\p\p}, {\bar{m}\phantom{|}\!}^{\p\p\p}\rangle. 
\end{equation}
Then they satisfy
\setlength\arraycolsep{1pt}
\[\begin{array}{rcl}
F_2&=&\displaystyle m^{\p\p}\cdot m^{\p\p}=\kappa^2,\\
F_3&=&\displaystyle m^{\p\p\p}\cdot m^{\p\p\p}=\kappa^4+{\kappa^{\p}}^2+\kappa^2\tau^2.
\end{array}\] \setlength\arraycolsep{5pt}
By derivating these equations we obtain a table of $\langle {\bar{m}}^{(i)}(0), {\bar{m}}^{(j)}(0) \rangle$ needed in this article (Table \ref{table_fifj}). 

\begin{table}[h]
\begin{center}
\par\noindent
\begin{tabular}{|c|c|c|c|c|c|}
\hline
& $\!\!\phantom{\stackrel{{a}}{0}}\!\! j=0\!\!$  &  $\!\!j=1\!\!$  &  $\!\!j=2\!\!$  & $\!\!j=3\!\!$ & $j=4$ \\[1mm] 
\hline
 $\!\!\phantom{\stackrel{{a}}{0}}\!\!i=0$ & $0$ & $0$ & $-1$ & $0$ & $F_2$ \\[1mm] 
\hline
 $\!\!\phantom{\stackrel{{a}}{0}}\!\!i=1$ &  & $1$ & $0$ & $-F_2$ & $\ast$ \\[1mm] 
\hline
 $\!\!\phantom{\stackrel{{a}}{0}}\!\!i=2$ &  &  & $F_2$ & $\ast$ &$\ast$ \\[1mm] 
\hline
 $\!\!\phantom{\stackrel{{a}}{0}}\!\!i=3$ &  &  &  & $F_3$ & $\ast$ \\[1mm]
\hline 
\end{tabular}
\end{center}
\caption{A table of $\langle {\bar{m}}^{(i)}(0), {\bar{m}}^{(j)}(0) \rangle$ \hspace{0.3cm}($F_2=\kappa^2$, $F_3=\kappa^4+{\kappa^{\p}}^2+\kappa^2\tau^2$)}
\label{table_fifj}
\end{table}

\smallskip
If $u<v<w$ then $\bar{m}(u), \bar{m}(v)$, and $\bar{m}(w)$ are linearly independent in $\mathbb R^5_1$, and therefore (\ref{S=Gr=Theta}) implies that $\gamma(u,v,w)$ is given by 
\begin{equation}\label{f_gamma(s,t,u)}
\gamma(u,v,w)=\frac{\bar{m}(u)\w \bar{m}(v)\w \bar{m}(w)}{\phantom{\hat{|}}\!\big\Vert\bar{m}(u)\w \bar{m}(v)\w \bar{m}(w)\big\Vert}.\end{equation}

\begin{lemma}\label{formula_circle_through_a_knot} 
Let $\Gamma(s)$ be an osculating circle to a curve $C=\{\bar{m}(s)\}$ in  $\EE^3$ which is parametrized by the arc-length $s$. 
Then $\Gamma(s)$ is given by 
\[\Gamma(s)=\EE^3\cap\spa\big\langle \bar{m}(s), {\bar{m}\phantom{|}\!}^{\p}(s), {\bar{m}\phantom{|}\!}^{\p\p}(s) \big\rangle.\]
The point $\gamma(s)$ in $\Theta(1,3)$ which corresponds to $\Gamma(s)$ is given by \setlength\arraycolsep{1pt}
\begin{equation}\label{f_gamma(s)}
\gamma(s)=\bar{m}(s)\w {\bar{m}\phantom{|}\!}^{\p}(s)\w {\bar{m}\phantom{|}\!}^{\p\p}(s). 
\end{equation}
\end{lemma}
\begin{proof} 
Consider the limit of $\gamma(u,v,w)$ as both $u$ and $w$ approach $v$. 
Taylor's expansion formula implies 
\[\begin{array}{l}
\bar{m}(s-\varDelta s)\w \bar{m}(s)\w \bar{m}(s+{\varDelta s}^{\p})\\[1mm]
=\displaystyle \left(\bar{m}(s)-\varDelta s{\bar{m}\phantom{|}\!}^{\p}(s)+\frac{(\varDelta s)^2}{2}{\bar{m}\phantom{|}\!}^{\p\p}(s)+O\left((\varDelta s)^3\right)\right)\w\bar{m}(s)\\[4mm]
\hspace{0.5cm}\displaystyle \w \left(\bar{m}(s)+{\varDelta s}^{\p}{\bar{m}\phantom{|}\!}^{\p}(s)+\frac{({\varDelta s}^{\p})^2}{2}{\bar{m}\phantom{|}\!}^{\p\p}(s)+O\left(({\varDelta s}^{\p})^3\right)\right)\\[6mm]
=\displaystyle \frac{\varDelta s{\varDelta s}^{\p}(\varDelta s+{\varDelta s}^{\p})}{2}\,
\bar{m}(s)\w {\bar{m}\phantom{|}\!}^{\p}(s)\w {\bar{m}\phantom{|}\!}^{\p\p}(s)+\mbox{higher order terms}\,.
\end{array}\]
It follows that the osculating circle is given by the vector $\bar{m}(s)\w {\bar{m}\phantom{|}\!}^{\p}(s)\w {\bar{m}\phantom{|}\!}^{\p\p}(s)$ multiplied by a positive number. 
By Formula (\ref{f_wedge-det_x}) and Table \ref{table_fifj} we have 
\setlength\arraycolsep{5pt}
\[\big\langle \bar{m}\w{\bar{m}\phantom{|}\!}^{\p}\w{\bar{m}\phantom{|}\!}^{\p\p}, \bar{m}\w{\bar{m}\phantom{|}\!}^{\p}\w{\bar{m}\phantom{|}\!}^{\p\p}\big\rangle
=-\left|\begin{array}{ccc}
0 & 0 & -1 \\
0 & 1 & 0 \\
-1 & 0 & m^{\p\p}\cdot m^{\p\p} 
\end{array}\right|
=1. \]
Therefore the osculating circle $\gamma(s)$ is given by 

\[\gamma(s)=\frac{\bar{m}(s)\w {\bar{m}\phantom{|}\!}^{\p}(s)\w {\bar{m}\phantom{|}\!}^{\p\p}(s)}{\big\Vert\bar{m}(s)\w {\bar{m}\phantom{|}\!}^{\p}(s)\w {\bar{m}\phantom{|}\!}^{\p\p}(s)\big\Vert}=\bar{m}(s)\w{\bar{m}\phantom{|}\!}^{\p}(s)\w{\bar{m}\phantom{|}\!}^{\p\p}(s).\]
\end{proof}

\begin{theorem}\label{thm_osc_circles} 
Let $\gamma$ be a curve in $\mathcal{S}(1,3)$ which corresponds to the set of osculating circles of a curve $C$ in $\EE^3$. 
Then the curve $\gamma$ is a lightlike curve. 
\end{theorem}
\begin{proof} 
The formula (\ref{f_gamma(s)}) implies
\begin{equation}\label{f_gamma_prime}
\gamma^{\p}(s)=\bar{m}(s)\w{\bar{m}\phantom{|}\!}^{\p}(s)\w{\bar{m}\phantom{|}\!}^{\p\p\p}(s).
\end{equation}
By Formula (\ref{f_wedge-det_x}) and Table \ref{table_fifj} we have \setlength\arraycolsep{5pt}
\[\big\langle \bar{m}\w{\bar{m}\phantom{|}\!}^{\p}\w{\bar{m}\phantom{|}\!}^{\p\p\p}, \bar{m}\w{\bar{m}\phantom{|}\!}^{\p}\w{\bar{m}\phantom{|}\!}^{\p\p\p}\big\rangle
= -\left|\begin{array}{ccc}
0 & 0 & 0 \\
0 & 1 & \ast \\
0 & \ast & \ast
\end{array}\right|=0,\] 
which implies that $|\gamma^{\p}(s)|\equiv 0$. 
This ends the proof of Theorem \ref{thm_osc_circles}.
\end{proof}

\subsection{Conformal arc-length via osculating circles}
Let $C=\{\bar m(t)\}$ $(0\le t\le T)$ be a curve in $\SS^3$ or $\EE^3$, and $\gamma$ a lightlike curve in $\mathcal{S}(1,3)$ which consists of the osculating circles of $C$. 
Let $I$ denote the domain interval of the map $\bar m$. 

\begin{theorem}\label{drho<->L(gamma'')ds}
By taking the pull-back to $I$, the conformal arc-length element $d\rho$ of $C$ 
is equal to $\sqrt[4]{12}$ times the $\frac12$ dimensional length element $d\rho_{L^{\frac12}(\gamma)}$ of $\gamma$: 
\begin{equation}\label{rho<->gamma''}
d\rho=\gamma^{\ast}\big(\sqrt[4]{12}\,d\rho_{L^{\frac12}(\gamma)}\big)
=\sqrt[4]{\LL(\Ppt{\gamma})}\,dt. 
\end{equation}
\end{theorem}

\begin{proof} 
We use the arc-length parameter $s$ of $C$. 
Let us abbreviate $\gamma^{(i)}(0)$ and ${\bar{m}\phantom{|}\!}^{(i)}(0)$ as $\gamma^{(i)}$ and ${\bar{m}\phantom{|}\!}^{(i)}$ in the proof. 
Since 
\[\gamma^{\p\p}= \bar{m}\w{\bar{m}\phantom{|}\!}^{\p}\w{\bar{m}\phantom{|}\!}^{(4)}+\bar{m}\w{\bar{m}\phantom{|}\!}^{\p\p}\w{\bar{m}\phantom{|}\!}^{\p\p\p} \,,\] 
Formula (\ref{f_wedge-det_x}) and Table \ref{table_fifj} imply 
\begin{eqnarray}
\LL\left(\gamma^{\p\p}\right)\!\!&\!\!=\!\!&\!\!\displaystyle\langle \gamma^{\p\p}, \gamma^{\p\p} \rangle \nonumber\\
\!\!&\!\!=\!\!&\!\!\displaystyle \big\langle\bar{m}\w{\bar{m}\phantom{|}\!}^{\p}\w{\bar{m}\phantom{|}\!}^{(4)}, \bar{m}\w{\bar{m}\phantom{|}\!}^{\p}\w{\bar{m}\phantom{|}\!}^{(4)}\big\rangle
+2\big\langle\bar{m}\w{\bar{m}\phantom{|}\!}^{\p}\w{\bar{m}\phantom{|}\!}^{(4)}, \bar{m}\w{\bar{m}\phantom{|}\!}^{\p\p}\w{\bar{m}\phantom{|}\!}^{\p\p\p}\big\rangle \nonumber\\
\!\!&\!\! \!\!&\!\!\displaystyle +\big\langle\bar{m}\w{\bar{m}\phantom{|}\!}^{\p\p}\w{\bar{m}\phantom{|}\!}^{\p\p\p}, \bar{m}\w{\bar{m}\phantom{|}\!}^{\p\p}\w{\bar{m}\phantom{|}\!}^{\p\p\p}\big\rangle \nonumber\\[2mm]
\!\!&\!\!=\!\!&\!\!\displaystyle 
-\left(\,
\left|\begin{array}{ccc}
0 & 0 & F_2\\
0 & 1 & \ast\\
F_2 & \ast & \ast
\end{array}\right|
+2\left|\begin{array}{ccc}
0 & -1 & 0 \\
0 & 0 & -F_2\\
F_2 & \ast & \ast
\end{array}\right|
+\left|\begin{array}{ccc}
0 & -1 & 0 \\
-1 & \ast & \ast \\
0 & \ast & F_3
\end{array}\right|
\,\right) \nonumber\\[2mm]
\!\!&\!\!=\!\!&\!\!F_3-F_2^{\,2} \label{f_gamma''_F_3-F_2^2} \\[1mm]
\!\!&\!\!=\!\!&\!\!{\kappa^{\p}}^2+\kappa^2\tau^2\,, \label{f_L(gamma'')_kappa_tau}
\end{eqnarray}
which implies (\ref{rho<->gamma''}) since $d\rho=\sqrt[4]{{\kappa^{\p}}^2+\kappa^2\tau^2\,}\,ds$ by (\ref{inf_conf_arc-length}) and $\gamma^{\ast}d\rho_{\!L^{\frac12}(\gamma)}=\sqrt[4]{\frac{|\LL({\gamma}^{\p\p})|}{12}} \,\,ds$ by (\ref{f_1/2-length_element}). 
\end{proof}

\begin{remark}
A pair of nearby osculating circles is a ``{\it timelike pair}'', i.e. $\gamma(s+\varDelta s)-\gamma(s)$ is timelike for $|\varDelta s|\ll 1$. 
\end{remark}

\begin{corollary}\label{thm_osc_circles2}
The conformal arc-length parameter of a point $\bar m(t)$ in a curve $C$ in $\SS^3$ or $\EE^3$ can be expressed using the $L^{\frac12}$-measure of an arc of the lightlike curve $\gamma$ in $\mathcal{S}(1,3)$ which consists of osculating circles to $C$:  
\begin{equation}\label{f_conf_arc-length<->osc_circle}
\rho(t)-\rho(0)=\sqrt[4]{12}\lim_{\max|t_{j+1}-t_j|\to+0}\sum_{i}\sqrt{\Vert\gamma(t_{i+1})-\gamma(t_i)\Vert}\,,
\end{equation}
where $0<t_1<\cdots<t_k=t$ is a subdivision of the interval $[0,t]$. 
\end{corollary}
\begin{corollary}
The conformal arc-length parameter can be characterized as the parameter that satisfies
$\langle \Ppt{\gamma}, \Ppt{\gamma}\,\rangle\equiv1.$ 
Namely, it is a parameter which makes the curvature of a curve of osculating circles being identically equal to $1$. 
\end{corollary}
\begin{definition}\label{def_F3-F2^2-form} \rm 
Let $C$ be a curve in $\RR^n$. 
Suppose it can be expressed by the arc-length as $C=\{m(s)\}$. 
Define the $1$-form $\omegasub{C}$ on $C$ by 
\[m^{\ast}\omegasub{C}=\sqrt[4]{m^{\p\p\p}\cdot m^{\p\p\p}-(m^{\p\p}\cdot m^{\p\p})^2 \,}\,ds,\]
where ${}^{\p}$ denotes the derivation with respect to the arc-length $s$. 
\end{definition}
\begin{corollary}\label{cor_inv_4v/F3-F2^2ds} 
As above the $1$-form $\omegasub{C}$ is invariant under M\"obius transformations. 
Namely, if $G$ is a M\"obius transformation then we have 
$G^{\ast}\omegasub{G(C)}=\omegasub{C}.$ 
When $n=3$ this $1$-form $\omegasub{C}$ is equal to the conformal arc-length element $d\rho$ of $C$. 
\end{corollary}
\begin{remark}\rm 
Theorem 1.1 of Liu (\cite{Liu}) implies that if $\bar C=\{\bar m(s)\}$ is a curve in the light cone such that ${\bar m}^{\p}$ is not parallel to $m$ then the $1$-form $\omegasub{\bar C}$ given by 
\[{\bar m}^{\ast}\omegasub{\bar C}=\sqrt[4]{\big\langle {\bar m}^{\p\p\p}, {\bar m}^{\p\p\p}\big\rangle-\big\langle {\bar m}^{\p\p}, {\bar m}^{\p\p}\big\rangle^2 \,}\,ds,\]
where $s$ is the arc-length, is invariant under any transformation of the form $T: \bar m(s)\mapsto e^{f(s)}\bar m(s)$  . 

Let $\iota:\RR^n\to\EE^n_0\subset\RR^{n+2}_1$ be the natural bijection (\ref{def_E^3_0}). 
Through this bijection a M\"obius transformation of $\RR^n$ can be expressed as $\iota^{-1}\circ T\circ A\circ \iota$ for some $A\in O(n+1,1)$ and some transformation $T$ of the previous form of the light cone. 
Since $\iota$ is an isometry and $A\in O(n+1,1)$ does not change $\omegasub{\bar C}$, Corollary \ref{cor_inv_4v/F3-F2^2ds} follows from Liu's result. 
\end{remark}
\begin{proof}
We prove it when $n=3$. 
Suppose $G\cdot C$ are expressed as $G(C)=\{m_2(t)\}$ with $t$ being the arc-length of $G(C)$. 
Let $\gamma_{m}$ and $\gamma_{m_2}$ be the curves in $\mathcal{S}(1,3)$ which consists of the osculating circles to $C$ and $G(C)$ respectively. 
They are lightlike curves by Theorem \ref{thm_osc_circles}. 
The formula (\ref{f_gamma''_F_3-F_2^2}) implies that 
\setlength\arraycolsep{1pt}
\[\begin{array}{rcl}
\sqrt[4]{m^{\p\p\p}\cdot m^{\p\p\p}-(m^{\p\p}\cdot m^{\p\p})^2 \,}\,ds&=&\displaystyle \sqrt[4]{\left\langle\gamma_{m}^{\>\,\p\p}, \gamma_{m}^{\>\,\p\p}\right\rangle}\,ds, \\ [2 mm]
\sqrt[4]{\Pppt m_2\cdot \Pppt m_2-(\Ppt m_2\cdot \Ppt m_2)^2 \,}\,dt&=&\displaystyle \sqrt[4]{\left\langle{\Ppt\gamma_{m_2}}, {\Ppt\gamma_{m_2}}\right\rangle}\,dt, 
\end{array}\]
where putting $\Pt{}$ above means taking the derivation with respect to $t$. 
Lemma \ref{prop_tilde_A} shows that the M\"obius transformation $G\in O(4,1)$ produces a pseudoorthogonal transformation $\widetilde G\in O(6,4)$, and $\gamma_{m_2}$ is given by $\gamma_{m_2}=\widetilde G\cdot \gamma_{m}$. 
Since $\langle{\Ppt\gamma_{m_2}}, \,{\Ppt\gamma_{m_2}}\rangle=\big\langle\widetilde G\cdot {\Ppt\gamma_{m}}, \,\widetilde G\cdot {\Ppt\gamma_{m}}\big\rangle=\langle{\Ppt\gamma_{m}}, \, {\Ppt\gamma_{m}}\rangle$, Corollary \ref{cor_inv_4v/F3-F2^2ds} is a consequence of Lemma \ref{lem_inv_form_null_curve}. 
\end{proof}
The above proof shows that the pull-back to $I$ of $\omegasub{C}$ is equal to that of the $\frac12$ dimensional length element $d\rho_{\!L^{\frac12}(\gamma)}$ of the curve of osculating circles, where $I$ is the domain interval of the curve $C$. 

Corollaries \ref{drho<->L(gamma'')ds} and \ref{cor_inv_4v/F3-F2^2ds} show that the conformal arc-length is in fact invariant under M\"obius transformations. 
\subsection{Characterization of vertices of space curves in terms of osculating circles}

\begin{definition}\label{def_vertex} \rm 
We say that a point is a {\em vertex} of a curve $C$ if the the osculating circle has the third order contact with $C$ at that point. 
\end{definition}

\begin{remark}\label{rem_vertex}\rm 
(1) The condition is same as the usual one, $k^{\p}=0$, for plane curves. 

(2) The notion of a vertex is conformally invariant. 

(3) By the strong transversality  principle (cf. \cite{arnold}), vertex-free curves are generical in the sense that they form an open dense subset in the $C^\infty$ topology.

(4) A different definition of vertex can be found in the literature. 
In some works vertices are points where the osculating sphere has the $4$th order contact (cf.\cite{uribe}). 
There are many examples of vertex-free curves in our sense having the $4$th order contact with its
osculating sphere. 
For instance, one can consider plane curves. 
\end{remark}

\medskip
Bouquet's formula says that a curve can be expressed, with a suitable coordinates around a point, as 
\setlength\arraycolsep{1pt}\begin{eqnarray}\label{Bouquet_formula}
\begin{array}{rcllr}
x&=&\displaystyle s&&-\frac{\kappa}6s^3+\cdots\\[1mm]
y&=&\displaystyle &\phantom{+}\frac{\kappa}2s^2&+\frac{\kappa^{\p}}6s^3+\cdots\\[1mm]
z&=&\displaystyle &&\phantom{+}\frac{\kappa\tau}6s^3+\cdots
\end{array}
\end{eqnarray}
where $s, \kappa$, and $\tau$ are the arc-length, curvature, and torsion respectively. 

It is easy to check the following lemma using Bouquet's formula: 
\begin{lemma}\label{lem_vertex_1} 
Let $C$ be a curve in $\mathbb R^3$. 
A point $x$ is a vertex of $C$ if and only if 
\[{\kappa^{\p}}^2+{\kappa}^2\tau^2=0\] 
holds at $x$, where $\kappa$ and $\tau$ are the curvature and torsion of the curve $C$ (the derivation is taken with respect to the arc-length parameter of $C$). 
\end{lemma}
Now a vertex can be characterized in terms of the curve of osculating circles. 
Suppose $\EE^3$ (or $\SS^3$) is an Euclidean (or respectively, a spherical) model in $\RR^5_1$. 
\begin{theorem}\label{thm_vertex-osc_circle}
Let $C=\{\bar m(t)\}$ be a curve in $\EE^3$ or $\SS^3$ and $\gamma(t)$ a point in $\Theta(1,3)$ which corresponds to the osculating circle to $C$ at $\bar m(t)$. 
Then the following three are equivalent conditions for a point $\bar m(t_0)$ to be a vertex of $C$: 
\begin{enumerate}
\item $\Pt{\gamma}(t_0)=\vect 0$. 
\item $\big\langle\Ppt{\gamma}(t_0), \Ppt{\gamma}(t_0)\big\rangle=0$, i.e. the $\frac12$ dimensional length element of $\gamma$ vanishes at $t_0$. 
\item $\dim\span{\Pt{\gamma}, \Ppt{\gamma}}<2$. 
\end{enumerate}
\end{theorem}

\begin{proof} 
(1) Suppose $m(s)$ is a curve in $\RR^3$ parametrized by the arc-length $s$, and $C=\{\bar m(s)\}$ is a corresponding curve in $\EE^3_0$ which is given by (\ref{m_bar_m}). 
Then we have \setlength\arraycolsep{1pt}
\[\begin{array}{rcccc}
\bar m&=\displaystyle \left(\phantom{\frac{m}4}\hspace{-0.4cm}\right. & \displaystyle 1+\frac{m\cdot m}4, &\displaystyle -1+\frac{m\cdot m}4, &\displaystyle \left. m\phantom{\frac{m}4}\hspace{-0.2cm}\right),\\[4mm]
{\bar m}^{\p}&=\displaystyle \left(\phantom{\frac{m}4}\hspace{-0.4cm}\right. & \displaystyle \frac{m\cdot m^{\p}}2, &\displaystyle \frac{m\cdot m^{\p}}2, &\displaystyle \left. m^{\p}\phantom{\frac{m}4}\hspace{-0.2cm}\right),\\[4mm]
{\bar m}^{\p\p}&=\displaystyle \left(\phantom{\frac{m}4}\hspace{-0.4cm}\right. & \displaystyle \frac{m\cdot m^{\p\p}+1}2, &\displaystyle \frac{m\cdot m^{\p\p}+1}2, &\displaystyle \left. m^{\p\p}\phantom{\frac{m}4}\hspace{-0.2cm}\right),\\[4mm]
{\bar m}^{\p\p\p}&=\displaystyle \left(\phantom{\frac{m}4}\hspace{-0.4cm}\right. & \displaystyle \frac{m\cdot m^{\p\p\p}}2, &\displaystyle \frac{m\cdot m^{\p\p\p}}2, &\displaystyle \left. m^{\p\p\p}\phantom{\frac{m}4}\hspace{-0.2cm}\right).
\end{array}\]
Suppose $t_0$ corresponds to $s=0$. 
We may assume, after a M\"obius transformation if necessary, that $\kappa(0)\ne0$. 
Lemma \ref{lem_vertex_1} implies that $\bar m(0)$ is a vertex if and only if ${\kappa^{\p}}^2+\kappa^2\tau^2=0$. 
Since 
\[m(s)=m(0)+\vect e_1s+\kappa\vect e_2\frac{s^2}2+\left(-\kappa^2\vect e_1+\kappa^{\p}\vect e_2+\kappa\tau\vect e_3\right)\frac{s^3}6+O(s^4),\]
$\bar m(0)$ is a vertex if and only if $m^{\p\p\p}(0)=-\kappa^2m^{\p}(0)$, which can occur if and only if ${\bar m}^{\p\p\p}(0)=-\kappa^2{\bar m}^{\p}(0)$. 
As $a{\bar m}(0)+b{\bar m}^{\p}(0)+c{\bar m}^{\p\p\p}(0)=\vect 0$ $(a,b,c\in\RR)$ always implies $a=0$, ${\bar m}(0)$ is a vertex if and only if  ${\bar m}(0), {\bar m}^{\p}(0)$, and ${\bar m}^{\p\p\p}(0)$ are linearly dependent, i.e. ${\bar m}(0)\w {\bar m}^{\p}(0)\w {\bar m}^{\p\p\p}(0)=\vect 0$. 
Since Formula (\ref{f_gamma(s)}) implies that the left hand side of the last equation is equal to $\gamma^{\p}(0)$, it means the equivalence of the condition (1) and $\bar m(t_0)$ being a vertex. 

\medskip
The condition (2) follows from Lemma \ref{lem_vertex_1} and Formula (\ref{f_L(gamma'')_kappa_tau}). 

\medskip
The condition (1) implies the condition (3). 

On the other hand, if $\dim\span{\Pt{\gamma}, \Ppt{\gamma}}<2$ and $\Pt{\gamma}(t_0)\ne\vect 0$ then $\Ppt{\gamma}(t_0)=c\Pt{\gamma}(t_0)$ for some $c\in\RR$, which implies the condition (2). 
\end{proof}

\subsection{Characterization of curves of osculating circles}
We identify the set $\mathcal{S}(1,3)$ of the oriented circles in $\SS^3$ with the submanifold $\Theta(1,3)$ of $\RR^{10}_4$ as before. 
We give a condition for a curve in $\mathcal{S}(1,3)$ to be a set of osculating circles to a curve in $\SS^3$. 
\begin{theorem}\label{thm_curve_osc_circles}
A curve $\gamma$ in $\mathcal{S}(1,3)$ is a set of the osculating circles to a vertex-free curve in $\SS^3$ with a non-vanishing velocity vector if and only if $\gamma$ satisfies the following three conditins: 
\begin{enumerate}
\item[{\rm (i)}] $\gamma$ is a lightlike curve, 
\item[{\rm (ii)}] $\Pt{\gamma}(t)$ satisfies the Pl\"ucker relations for all $t$, 
\item[{\rm (iii)}] $\dim\spanbig{\Pt{\gamma}(t), \Ppt{\gamma}(t)}=2$. 
\end{enumerate}

Without the two conditions, $C$ being vertex-free and $\dim\span{\Pt{\gamma}, \Ppt{\gamma}}=2$, the ``if'' part of the above statement may fail although the ``only if'' part still holds. 

\end{theorem}
\begin{proof} 
We prove the first statement first. 

(1-1) ``Only if\,'' part. 
The condition (i) follows from Theorem \ref{thm_osc_circles}, the condition (ii) from Formula (\ref{f_gamma_prime}), as it means that $\gamma^{\p}$ is a pure $3$-vector and therefore satisfies the Pl\"ucker relations, and the condition (iii) follows from Theorem \ref{thm_vertex-osc_circle}. 

\medskip
(1-2) ``If\,'' part. 
Lemma \ref{lemm_time_cond} implies that $\Pt{\gamma}$ corresponds to an isotropic $3$-space tangent to the light cone. 
Therefore it defines a point in $\EE^3_0$, which we denote by $\bar m(t)$. 
Since $\Pt{\gamma}(t)$ and $\Ppt{\gamma}(t)$ are linearly independent, $\Pt{\bar m}(t)$ does not vanish. 
So, we may assume that $t$ is equal to the arc-length $s$ of $\bar m$. 
Define a map $m$ to $\RR^3$ by (\ref{m_bar_m}). 
We may assume, after a M\"obius transformation if necessary, that $m^{\p\p}(s_0)\ne\vect 0$. 

Put 
\[\bar{n}=\left(\frac{m\cdot n}2, \, \frac{m\cdot n}2, \, n\right)\in\RR^5_1 \hspace{0.3cm}\textrm{where}\hspace{0.3cm}n=\frac{m^{\p}\times m^{\p\p}}{\Vert m^{\p}\times m^{\p\p} \Vert}\in\RR^3.\]
Then $\bar m, \bar m^{\p}, \bar m^{\p\p}, \bar n$, and $n_1=(1,1,0,0,0)$ are linearly independent in $\RR^5_1$. 
Therefore, 10 pure $3$-vectors obtained as exterior products of three of $\bar m, \bar m^{\p}, \bar m^{\p\p}, \bar n$, and $n_1$ are linearly independent in $\RR^{10}_4$. 
We need formulae which express $\gamma(s)$ and $\gamma^{\p}(s)$ as linear combinations of them. 

\smallskip
For this purpose, we have to show that the point ${m}(s)$ belongs to the circle $\Ga(s)$ which corresponds to $\gamma(s)$. 

This can be proven as follows. 
We may assume, without loss of generality, that $\Ga(s)$ can be obtained as the intersection of our model of $\ss^3$ or $\EE ^3$ in the light cone and $\span{\vect e_0, \vect e_1, \vect e_2}$. 
Then, a computation shows that if $\gamma^{\p}$ satisfies the two conditions (i) and (ii) then there is a {\em lightlike pencil} $\mathcal{P}$ through $\gamma(s)$ (that is a $1$-parameter family of the oriented circles that are contained in a $2$-sphere which contains $\Ga(s)$ and that are all tangent to $\Ga(s)$ at a point $Q$) such that $\gamma^{\p}(s)$ is equal to a tangent vector to $\mathcal{P}$ at $\gamma(s)$. 
It means that the tangent point $Q$ above mentioned can be obtained as the intersection of $\SS^3$ and the isotropic $3$-space that corresponds to $\gamma^{\p}(s)$, namely $Q$ is equal to $m(s)$, which implies $m(s)\in \Ga(s)$ (see \cite{La-OH2} for pencils). 

\smallskip
Since $m(s)\in \Ga(s)$, $\gamma(s)$ can be expressed as $\gamma(s)=\alpha(s)\bar m(s)\w u_0(s)\w v_0(s)$ for some $\alpha(s)\in\mathbb R$ and $u_0(s), v_0(s)\in\RR^5_1$. 
Therefore $\gamma(s)$ can be expressed as the linear combination of the following 6 pure $3$-vectors; 
\[\bar m \w \bar m^{\p} \w \bar m^{\p\p}, \bar m \w \bar m^{\p} \w \bar n, \bar m \w \bar m^{\p} \w n_1, \bar m \w \bar m^{\p\p} \w \bar n, \bar m \w \bar m^{\p\p} \w n_1, \>\mbox{and} \>\> \bar m \w \bar n \w n_1.\] 
If one of the three coefficients of the latter three does not vanish, there is a non-zero coefficient of $\bar m^{\p} \w u_1 \w v_1$ $(\{u_1, v_1\}\subset\{\bar m^{\p\p}, \bar n, n_1\})$ of $\gamma^{\p}(s)$. 
On the other hand, by the assumption (ii) of the theorem and the definition of $\bar m$, $\gamma^{\p}(s)$ can be expressed in the form 
\[\gamma^{\p}(s)=b(s)\,\bar m(s)\w u(s)\w v(s) \hspace{0.3cm}\left(b(s)\in\RR,\,u(s), v(s)\in (\span{\bar m(s)})^{\perp}\subset \RR^5_1\,\right),\]
which is a contradiction. 

Therefore, $\Ga(s)$ is tangent to $C$ at $\bar m(s)$, i.e. $\gamma(s)$ is of the form 
\[\gamma(s)=\xi(s)\,\bar m \w \bar m^{\p} \w \bar m^{\p\p}
+\eta(s)\,\bar m \w \bar m^{\p} \w \bar n +\zeta(s)\,\bar m \w \bar m^{\p} \w n_1.\]
By patient computation we get 
\[\left\langle\gamma^{\p}, \gamma^{\p}\right\rangle=\left(\xi^{\p}+2\zeta^{\p}\right)^2+\eta^2+4\kappa^2\zeta^2,\]
where $\kappa$ is the curvature of $C$, $\Vert m^{\p\p}\Vert$. 
Since $\langle\gamma^{\p}(s), \gamma^{\p}(s)\rangle=0$ and $\kappa\ne0$ by our assumption, we have $\eta(s)=\zeta(s)=0$. 
Then, $\langle\gamma(s), \gamma(s)\rangle=1$ implies $\xi(s)=1$, which completes the proof of the first statement. 

\medskip
(2) The second satement of the theorem can be verified in the same way as in Theorem \ref{lem_lightlike_curve=osc.cirlces}. 
\end{proof}

{The authors thank Martin Guest and Fran Burstall for informing the second author of the following condition due to Burstall. } 

\begin{proposition}{\rm (Burstall condition)} 
Suppose $\gamma\in\mathcal{S}(1,3)$ corresponds to a timelike $3$-space $\Pi$ of $\RR^5_1$. 
Recall that $T_{\gamma}\mathcal{S}(1,3)$ can be identified with $\textrm{Hom}(\Pi, \Pi^{\perp})$. 
Suppose an element $A$ in $\textrm{Hom}(\Pi, \Pi^{\perp})$ that corresponds to $\Pt{\gamma}$ can be expressed as \setlength\arraycolsep{2pt}
$A=\displaystyle \left(\begin{array}{ccc} a&b&c\\d&e&f \end{array}\right)=\left(\begin{array}{c} \vect a_1\\ \vect a_2 \end{array}\right)$ with respect to orthonormal bases of $\Pi$ and $\Pi^{\perp}$. 
Then the conditions {\rm (i)} and {\rm (ii)} of {\rm Theorem \ref{thm_curve_osc_circles}} are equivalent to the condition $\langle A, {}^t\!A\rangle=O$, namely, 
\[\langle\vect a_1, \vect a_1 \rangle=\langle\vect a_2, \vect a_2 \rangle=\langle\vect a_1, \vect a_2 \rangle=0,\]
where $\vect a_1$ and $\vect a_2$ are considered as vectors in the Minkowski space $\RR^3_1$. 
\end{proposition}

\begin{proof} 
We may assume without loss of generality that $\Pi(t_0)=\span{\vect e_0, \vect e_1, \vect e_2}$, i.e. $\gamma(t_0)=(0,\,\cdots,\,0,1)$. 
Then $\Pi^{\perp}(t_0)=\span{\vect e_3, \vect e_4}$. 
If we use $\{\vect e_0, \vect e_1, \vect e_2\}$ and $\{\vect e_3, \vect e_4\}$ as bases of $\Pi(t_0)$ and $\Pi^{\perp}(t_0)$ respectively, 
the Pl\"ucker coordinates of $\Pt{\gamma}(t_0)$ are given by  

\setlength\arraycolsep{2pt}
\[\begin{array}{c}
(p_{234}, \, p_{134}, \, p_{124}, \, p_{123}\, ; \, p_{034}, \, p_{024}, \, p_{023}, \, p_{014}, p\, _{013}, \, p_{012})\big(\Pt{\gamma}(t_0)\big)\\[1mm]
=(\,0,\,0,d,a;\,0,{-e},{-b},f,c,\,0\,). \end{array}
\]
Then the conditions (i) and (ii) of Theorem \ref{thm_curve_osc_circles} are given by 
\[\begin{array}{c}
-a^2-d^2+b^2+c^2+e^2+f^2=0,\\[2mm]
\left\{
\begin{array}{rcl}
bd-ae&=&0,\\
-cd+af&=&0,\\
ce-bf&=&0,
\end{array}\right.
\end{array}\]
which are equivalent to 
\[\begin{array}{c}
-a^2+b^2+c^2=-d^2+e^2+f^2=0,\\[1mm]
(a,b,c)/\!/ (d,e,f), 
\end{array}\]
which, in turn, are equivalent to (see Lemma 9.1.1 (2) of \cite{OH2}) 
\[\langle\vect a_1, \vect a_1 \rangle=\langle\vect a_2, \vect a_2 \rangle=\langle\vect a_1, \vect a_2 \rangle=0.\]
\end{proof}
\subsection{Conformal arc-length and conformal angles}

\begin{definition}\label{def_conf_angle}\rm (Doyle and Schramm) 
Let $x$ and $y$ be a pair of distinct points on a curve $C$. 
Let $\theta_{C}(x,y)$ $(0\le\theta_{C}(x,y)\le\pi)$ be the angle between $\Ga(x,x,y)$ and $\Ga(x,y,y)$. 
We call it the {\em conformal angle} between $x$ and $y$. 
\begin{figure}
\begin{center}
\includegraphics[width=8cm]{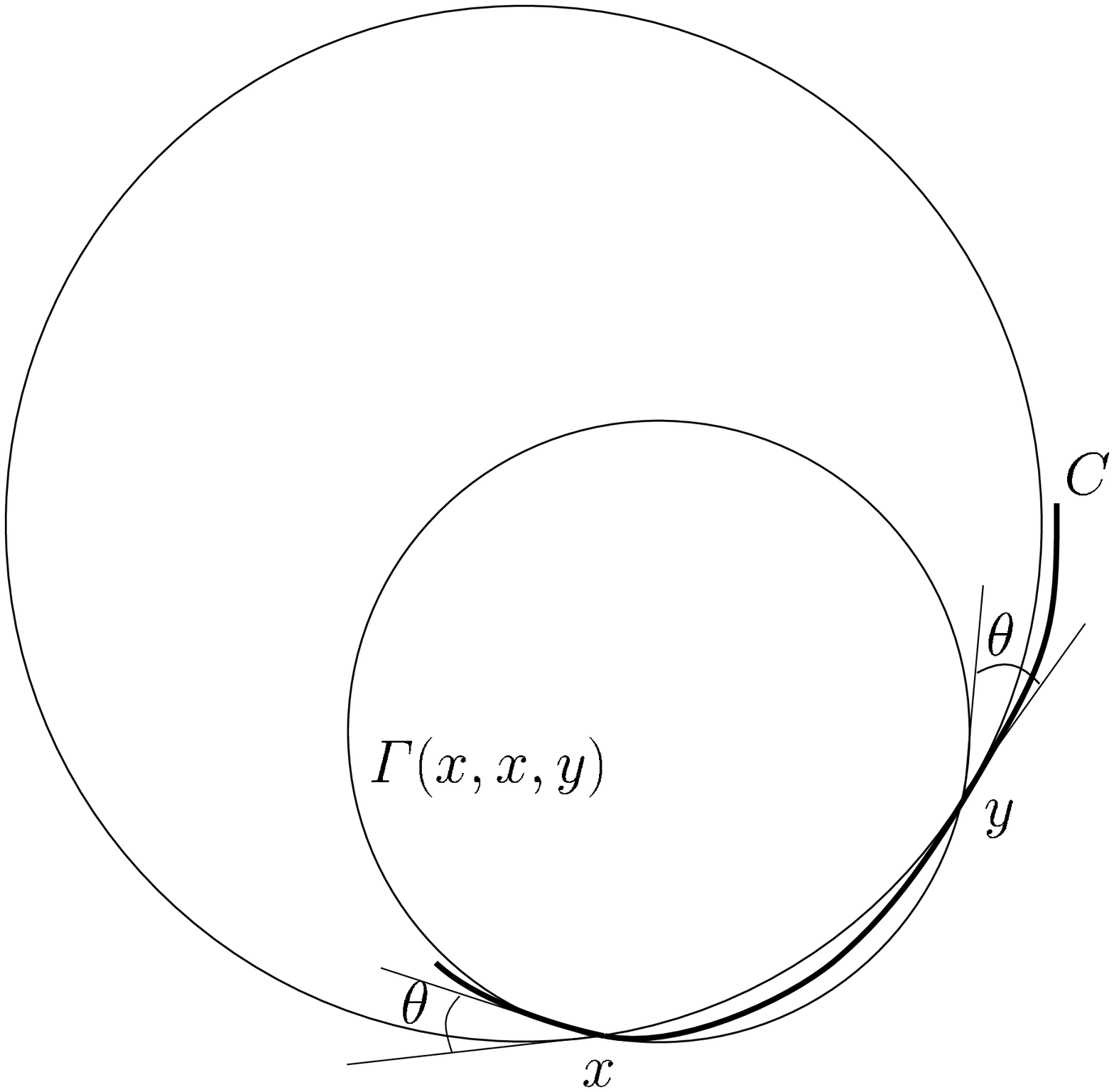} 
\caption{The conformal angle $\theta_C(x,y)$ \label{theta_C}. }
\end{center}
\end{figure}
\end{definition}

We note that the conformal angle is conformally invariant because it can be defined by angles, circles, and tangency, which are preserved by M\"obius transformations. 

Applying Bouquet's formula (\ref{Bouquet_formula}) to $\sin\theta_C$ we have 

\begin{lemma}{\rm (\cite{La-OH1})}\label{lemma_conf_angle} 
Let $s, \kappa, \tau$ be the arc-length, curvature, and torsion of $C$ respectively. 
Then the conformal angle satisfies 
\begin{eqnarray}\label{conf_angle_diagonal}
\theta _C(x,y)=\frac{\sqrt{{\kappa^{\p}}^2+\kappa^2\tau^2\,}}6\,|x-y|^2+O(|x-y|^3). 
\end{eqnarray}
\end{lemma}
The formula (\ref{inf_conf_arc-length}) of the conformal arc-length implies that the conformal arc-length can be interpreted in terms of the conformal angle as follows. 
\begin{proposition}\label{prop_cal<->cangle}
The conformal arc-length $\rho$ satisfies 
\[\frac{d\rho}{ds}(s)=\lim_{\varDelta s\to0}\frac{\sqrt{6\,\theta_{{C}}\big(m(s), m(s+\varDelta s)\big)}}{\varDelta s}.\]
\end{proposition}

\begin{remark}\rm 
$\>$Put 
\[\Pi=\textsl{Span}\langle \gamma(s,s,s+\varDelta s), \gamma(s,s+\varDelta s,s+\varDelta s)\rangle.\]
It is a spacelike $2$-plane of $\RR^{10}_4$. 
We claim that the intersection of $\Pi$ and $\Theta(1,3)$ is $1$ dimensional, moreover it is a circle which is a geodesic of $\Theta(1,3)$. 
The conformal angle is equal to the distance in $\mathcal{S}(1,3)$ (i.e. the shorter arc-length along this geodesic) between $\gamma(s,s,s+\varDelta s)$ and $\gamma(s,s+\varDelta s,s+\varDelta s)$. 

The above statements can be understood in the timelike $4$-space 
$$W=\span{m(s), m^{\p}(s), m(s+\varDelta s), m^{\p}(s+\varDelta s)}\subset\RR^5_1$$
 as all the events take place in $W$. 
Remark that $W$ intersects $\SS^3$ in a ``bitangent'' sphere $\Si(m, m+\varDelta s)$ that contains the two tangent circles to the curve $C$, $\Ga(m(s), m(s), m(s+\varDelta s))$ and $\Ga(m(s), m(s+\varDelta s), m(s+\varDelta s))$. 
Therefore, the intersection of $\Theta(1,3)$ and the set of oriented circles in the sphere $\Si(m, m+\varDelta s)=W\cap \SS^3$ can be isometrically identified with $3$-dimensional de Sitter space $\La^3$ which consists of oriented circles in $\Si(m, m+\varDelta s)$. 
Through this identification, $\Pi\cap\Theta(1,3)$ can be identified with the intersection of a spacelike $2$-plane in $\RR^5_1$ that corresponds to $\Pi$ and $\La^3$. 
It consists of the oriented circles in $\Si(m, m+\varDelta s)$ that pass through both $m(s)$ and $m(s+\varDelta s)$. 
It is a circle which is a geodesic of $\La^3$. 
The distance between any pair of points on this circle is equal to the angle between a pair of corresponding circles (see Theorem 9.5.1 of \cite{OH2} or \cite{HJ}). 
\end{remark}

\medskip
\begin{remark}\rm 
We have two kinds of infinitesimal interpretation of the conformal arc-length as the distance between a pair of nearby circles: 
one by osculating circles $\Ga(x,x,x)$ and $\Ga(y,y,y)$ (Theorem \ref{thm_osc_circles2}) and the other by tangent circles $\Ga(x,x,y)$ and $\Ga(x,y,y)$ (Proposition \ref{prop_cal<->cangle}). 
We remark that the former is a ``timelike'' pair and the latter a ``spacelike'' pair. 
\end{remark}

\section{Integral geometric viewpoint}\label{sec_last_section}
\subsection{Information from  two nearby osculating circles}

\begin{definition}\rm (\cite{La-OH1}) 
Let $P_1, P_2, P_3$, and $P_4$ be points on an oriented sphere $\Si$. 
They can be considered as complex numbers through an orientation preserving stereographic projection from $\Si$ to $\CC\cup\{\infty\}$. 
The cross ratio of $P_1, P_2, P_3$, and $P_4$ can be defined by that of the corresponding four complex numbers. 
We remark that it does not depend on the stereographic projection that is used to define it. 
\end{definition}
Let $C$ be a curve in $\ss^3$ or $\RR^3$, and $\mathcal S$ be a sphere which intersects $C$ orthogonally at a point $m(t_1)$ that passes through a nearby point $m(t_2)$. 
Let ${\cal O}_{m(t_i)}$ $(i=1,2)$ be the osculating circle to $C$ at $m(t_i)$, and $y(t_i)$ the intersection point of ${\cal O}_{m(t_i)}$ and $\mathcal S$ so that ${\cal O}_{m(t_i)}\cap\mathcal S=\{m(t_i), y(t_i)\}$. 
\begin{figure}[ht] 
\begin{center}
\psfrag{mt}{ $m(t_1)$}
\includegraphics[width=10cm]{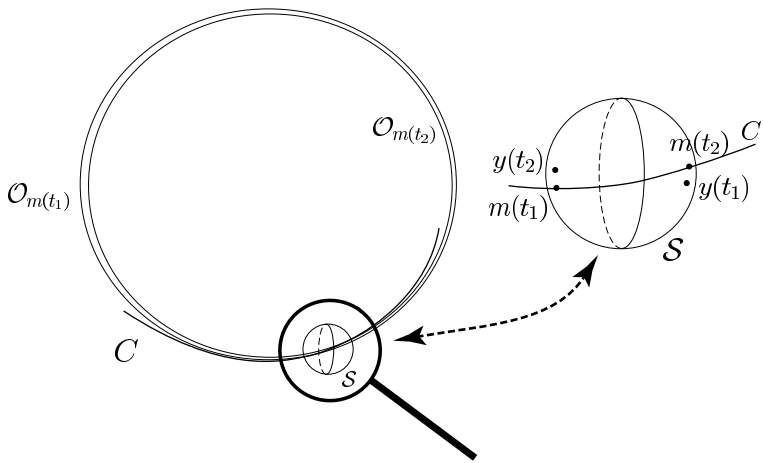}
\label{trace_on_S}
\caption{The intersection points of two nearby osculating circles and an almost orthogonal small sphere $\mathcal S$. In Proposition \ref{prop_cross-ratio_conf_arc-length}, $t_1$ is $t$ and $t_2$ is $t+dt$. }
\end{center}
\end{figure}
\begin{proposition}\label{prop_cross-ratio_conf_arc-length}
The infinitesimal cross ratio 
\[cross(m(t),y(t+dt);m(t+dt),y(t))=\frac{m(t)-y(t+dt)}{m(t)-y(t)}:\frac{m(t+dt)-y(t+dt)}{m(t+dt)-y(t)}\]
is real. 
The fourth root of its absolute value is equal to $\frac1{\sqrt6}$ times the pull-back of the conformal arc-length element of the curve $C$.
\end{proposition}
%
%
%

We remark that the infinitesimal cross ratio mentioned above is different from that defined in \cite{La-OH1} and studied in \cite{OH4}. 

\begin{proof}
Since both the cross ratio and the conformal arc-length element are invariant under M\"obius transformations, we may assume that $m(t_1)$ is the origin and that the curve $C$ in $\RR^3$ is given by the normal form (see \cite{CSW}) 
\begin{equation}\label{normal_form}
\begin{array}{l}
\displaystyle y=\frac{x^3}{3!} + O(x^5), \\[2mm]
z=0+O(x^4).\end{array}
\end{equation}
The normal form above can be obtained from Bouquet's formula (\ref{Bouquet_formula}) by putting $k=0$ and $k^{\p}=1$ at the origin. 

Let $h$ be the diameter of $\mathcal S$.  
Since the conformal arc-lenth element is given by $\sqrt[4]{(k')^2 + k^2 \tau^2}ds$, where $s$ is the arc-length of the curve, it is enough to show that the cross ratio is of the form $\frac{h^4}{36}+O(h^5)$. 

Since both $\Vert m(t)-y(t+dt)\Vert$ and $\Vert m(t+dt)-y(t)\Vert$ are of order $h^3$, we can neglect $O(h^4)$ terms. 
Therefore, the $z$-coordinate in the normal form (\ref{normal_form}) can be neglected. 
Then the four points $m(t), m(t+dt), y(t)$, and $y(t+dt)$ are on a circle that is the intersection of $\mathcal S$ and the $xy$-plane, which implies that the infinitesimal cross ratio is real. 

Let us now consider in the $xy$-plane. 
The coordinates of $m(t+dt)$ are given by $(h, \frac{h^3}{6})$ up to $O(h^4)$. 
A computation shows that the $x$-coordinate of the center of the osculating circle to the curve at $(h, \frac{h^3}{6})$ is $\frac{h}2+O(h^5)$, which implies that the intersection point $y(t+dt)$ of this osculating circle and $\mathcal S$ is equal to $(0, \frac{h^3}{6})$ up to $O(h^4)$. 
It follows that the absolute value of the infinitesimal cross ratio is given by 
\[\frac{\Vert m(t)-y(t+dt)\Vert }{\Vert m(t)-y(t)\Vert }\times \frac{\Vert m(t+dt)-y(t)\Vert }{\Vert m(t+dt)-y(t+dt)\Vert }=\frac{h^4}{36}+O(h^5),\]
which completes the proof. 
\end{proof}

\subsection{Integral geometric interpretation of the conformal arc-length element}\label{light_like_again}%
%
%
In the previous section we express the conformal arc-length as the $L^{\frac12}$-measure of a lightlike curve in the space of circles. 
The goal of this subsection is Theorem \ref{average_light}, where the conformal arc-length of a space curve $C$ is expressed as the average of the $L^{\frac12}$-measures of $1$ parameter family of lightlike curves in the space of spheres $\La^4$ that can be obtained from the curve of osculating spheres to $C$. 

The statement of Theorem \ref{average_light} is analogous to the following statement in the sense that something can be expressed as the average of $1$ parameter family of other quantities. 
\begin{proposition}\label{proj_eucl_curve}
The curvature of a curve $C\subset \RR^3$ at a point $m$ is proportional to the average of the curvatures  at $m$ of the plane curves obtained as the orthogonal projections of $C$ on the planes containing the tangent line $T_m C$ to $C$ at $m$ (the proportionality coefficient is $\pi$).  
\end{proposition}
\begin{proof}
We only need to project an osculating circle $\mathcal{O}_m$ to $C$ at $m$ on the planes containing $T_m C$ and observe the curvatures of these projections at the point $m$.
\end{proof}

\subsubsection{Preliminary lemmas}
Theorem \ref{average_light} has two kinds of proofs, a geometric one using pencils of spheres, pencils of circles and cross ratios, and an algebraic one using an ``anti-isometry'' $F$ (where $F$ being an anti-isometry means $\langle F(u), F(v)\rangle=-\langle u, v\rangle$ for any $u, v$) between two Grassmann manifolds. 
We start with preliminary lemmas which are needed for the geometric proof. 
\begin{definition} \rm 
The {\em Lorentz distance} between a pair of spheres $\Si_1$ and $\Si_2$ is the length of the geodesic $\gamma$ joining the two correspoinding points $\sigma_1$ and $\sigma_2$ in $\Lambda$: $\displaystyle d_{\LL}(\sigma_1,\sigma_2)=\int \Vert \Pt{\gamma}\Vert dt$. 
\end{definition}

Since a geodesic in $\La$ can be obtained as the intersection of $\Lambda$ with a $2$ dimensional vector subspace of $\RR^5_1$, the geodesic $\gamma$ joining $\sigma_1$ and $\sigma_2$ is a subarc of $\La\cap\span{\sigma_1, \sigma_2}$. 

The Lorentz distance between a pair of spheres can be expressed by the Lorentz distance between their intersections with an orthogonal sphere or an orthogonal circle. 
\begin{lemma} 
Let $d_{\LL}(\sigma_1,\sigma_2)$ be the Lorentz distance between a pair of spheres $\Si_1$ and $\Si_2$. 
\begin{enumerate}
    \item Let $S$ be a sphere orthogonal to $\Sigma_1$ and $\Sigma_2$. 
Then $d_{\LL}(\sigma_1,\sigma_2)$ is equal to the Lorentz distance $d_{\LL}(\gamma_1,\gamma_2)$ between the two circles $\Gamma_1= \Sigma_1 \cap S$ and $\Gamma_2=\Sigma_2 \cap S$. 
  \item Let ${\mathcal G}$ be a circle orthogonal to $\Sigma_1$ and $\Sigma_2$. 
Then $d_{\LL}(\sigma_1,\sigma_2)$ is equal to the Lorentz distance $d_{\LL}({\mathcal P}_1,{\mathcal P}_2)$ between the two $0$-spheres ${\mathcal P}_1= \Sigma_1 \cap {\mathcal G}$ and ${\mathcal P}_2=\Sigma_2 \cap {\mathcal G}$.
 \end{enumerate}
  \end{lemma}
\begin{figure}[ht] 
\begin{center}
\includegraphics[width=12cm]{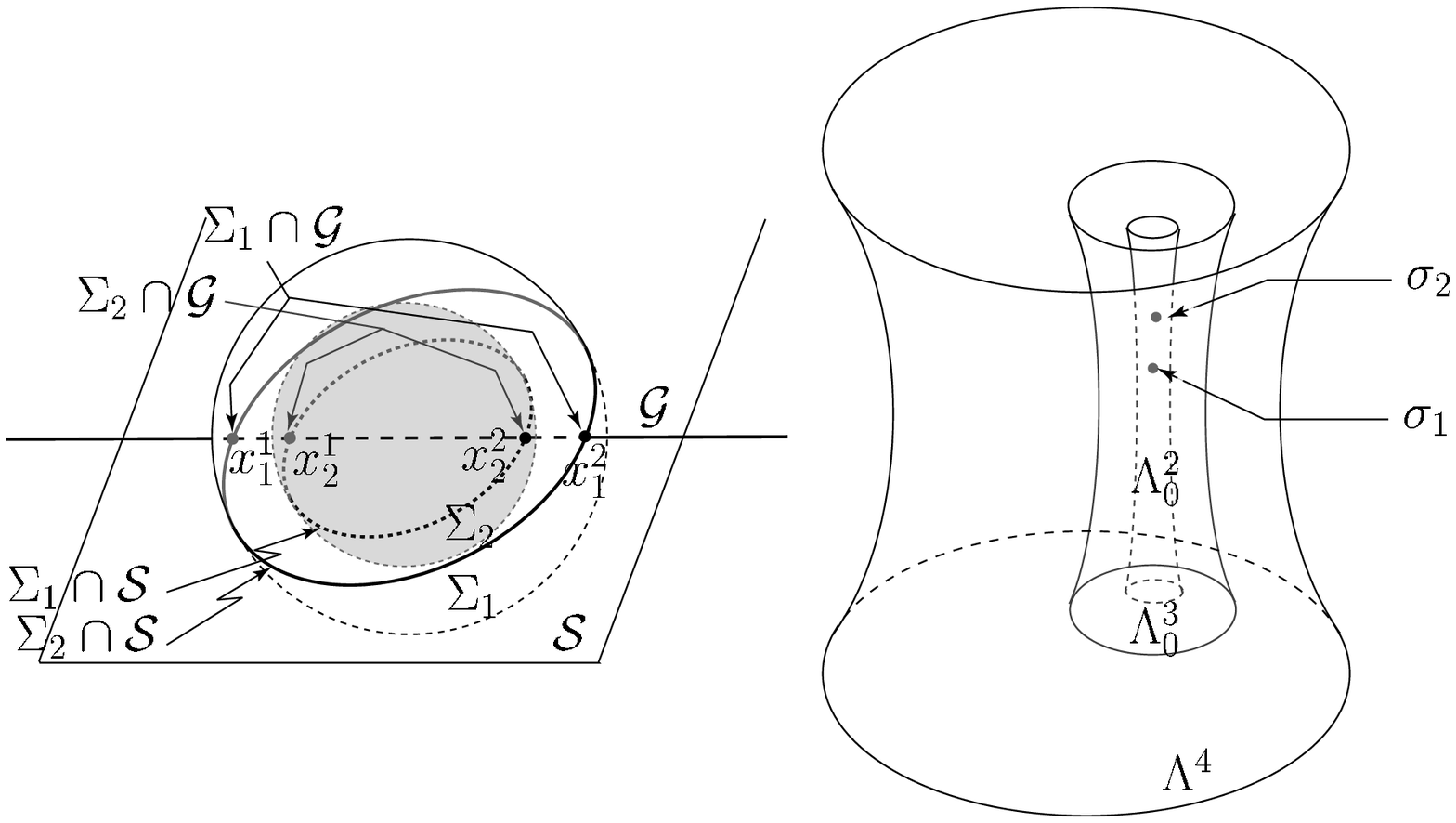}
\caption{Intersection and Lorentz distance \label{cross_lorentz}}
\end{center}
\end{figure}
\begin{proof} 
Let $\sigma_1$ and $\sigma_2$ be the two points in de Sitter space $\La^4$ which correspond to the two spheres $\Sigma_1$ and $\Sigma_2$. 

(1) can be proven by considering the two points $\sigma_1$ and $\sigma_2$ in de Sitter spaces $\Lambda^3_0 \subset  \Lambda^4$, where $\Lambda^3_0$ is the set of the oriented spheres orthogonal to the sphere $S$. 

To be more precise, let $\Pi$ be a $4$-dimensional vector subspace of $\RR^5_1$ so that $S$ is the intersection of $\Pi$ and $\SS^3$ or $\EE^3$, namely $\Pi=\spa{S}$, and $\sigma_0$ be the point in $\La^4$ which corresponds to $S$. 
Then $\La^3_0=\La^4\cap(\span{\sigma_0})^{\perp}=\La^4\cap\Pi$. 
Since the geodesic $\gamma$ in $\La^4$ joining $\sigma_1$ and $\sigma_2$ is a subarc of $\La^4\cap\span{\sigma_1, \sigma_2}$, it is contained in $\La^3_0$, which proves (1). 

\smallskip
(2) can be proven similarly, using de Sitter space $\Lambda^2_0$ which is the set of the oriented spheres orthogonal to the circle $\mathcal{G}$ (Figure \ref{cross_lorentz}). 

It can also be proven by the composition of (1) and a one dimensional lower analogue of (1). 
\end{proof}
\begin{remark}
The Lorentz distance between $\Sigma_1$ and $\Sigma_2$ and the cross ratio  of the four intersection points of ${\mathcal G}$ and $\Sigma_1 \cup \Sigma_2$ are related by a diffeomorphism. When the circle ${\mathcal G}$  is a line the cross ratio is given by the formula: 
$$cross(x_1^1,x_2^1;x_1^2,x_2^2)=\frac{x_1^1-x_2^1}{x_1^1-x_2^2}:\frac{x_1^2-x_2^1}{x_1^2-x_2^2}.$$
When the four points are on a circle in a complex plane, the four points $x_i^j$ should be considered as complex numbers, but the cross ratio is real as the points are on a circle.

 For example if we name the point as in Figure (\ref{cross_lorentz}), we have:
$$|cross(x_1^1,x_2^1;x_1^2,x_2^2)|=\left(\frac{e^{\ell} -1}{e^{\ell}+1}\right)^2 ,$$
\noindent where $\ell$ is the Lorentz distance $d_{\LL}(\sigma_1,\sigma_2)$.

In particular, when the Lorentz distance between $\Sigma_1$ and $\Sigma_2$ is small, we have:
 $$|cross(x_1^1,x_2^1;x_1^2,x_2^2)|\simeq \left(\frac{\ell}{2}\right)^2.$$
\end{remark}
We need another lemma to compare the Lorentz distance between a pair of spheres to the Lorentz distance between their intersections with a sphere or a circle almost orthogonal to them. 
\begin{lemma}\label{cross_almost} 
Let $\Si_1$ and $\Si_2$ be a pair of spheres. 
Let $cross(x_1^1,x_2^1;x_1^2,x_2^2)$ be the cross ratio of the intersection of $\Sigma_1$ and $\Sigma_2$ with a circle ${\mathcal G}$ which intersects $\Sigma_1$ and $\Sigma_2$ in the right angles. 
\begin{enumerate}
\item Let $cross(y_1^1,y_2^1;y_1^2,y_2^2)$ be the cross ratio of the intersection of $\Sigma_1$ and $\Sigma_2$ with a circle $\bar{\mathcal G}$ that makes angles $\theta_1$ with $\Sigma_1$ and $\theta_2$ with $\Sigma_2$. 
Then the quotient of $cross(x_1^1,x_2^1;x_1^2,x_2^2)$ and $cross(y_1^1,y_2^1;y_1^2,y_2^2)$ is in an interval $[1-\delta_1, 1+\delta_1]$, where $\delta_1$ is a function of $\theta_1$ and $\theta_2$ which goes to $0$ when $\theta_1$ and $\theta_2$ go to $\pm \pi /2$.

\item 
Suppose that a sphere $S$ makes angles $\theta_1$ with $\Sigma_1$ and $\theta_2$ with $\Sigma_2$. 
Let $cross(z_1^1,z_2^1;z_1^2,z_2^2)$ be the cross ratio of the intersection of $\Sigma_1 \cap S$ and $\Sigma_2 \cap S$ with a common orthogonal circle $\Gamma \subset S$. 
Then the quotient of $cross(x_1^1,x_2^1;x_1^2,x_2^2)$ and $cross(z_1^1,z_2^1;z_1^2,z_2^2)$ is in an interval of the form $[1-\delta_2, 1+\delta_2]$, where $\delta_2$ is  a function of $\theta_1$ and $\theta_2$ which goes to $0$ when $\theta_1$ and $\theta_2$ go to $\pm \pi /2$.
\end{enumerate}
\end{lemma}

\begin{proof}
Let us fix the two spheres $\Sigma_1$ and $\Sigma_2$.

Let us consider the circle $\tilde{\mathcal G}$  orthogonal to $\Sigma_1$ at $y^1_1$ and $y^1_2$. 
It intersects ${\mathcal G}$ at $y^1_1$ and $y^1_2$ making a small angle $\pi/2-\theta_1$.
It also intersects $\Sigma_2$ at points $z_1$ and $z_2$ close to $y^2_1$ and $y^2_2$ and with an angle close to $\pi/2$. 
All the corresponding arcs $a_i, b_i$ on both circles have a ratio satisfying $1-\delta <a_i/b_i <1+\delta$, which implies (1).  

The second statement (2) comes from the fact that the circle $\gamma$ is almost orthogonal to the spheres $\Sigma_1$ and $\Sigma_2$. 
\end{proof}
%
%
\subsubsection{Conformal arc-length as the average of $L^{\frac12}$-measure of lightlike curves in $\La^4$ which are associated to a curve in $\ss^3$ or $\RR^3$}

Let $C=\{c(t)\}$ be vertex-free curve in $\SS^3$ or $\EE^3$. 
Let $\Ga=\{\sigma(t)\}$ be a curve in $\La^4$ which is the set of osculating spheres of $C$. 
Let $\ts$ ($\ts$ is in some interval $I$) be the arc-length of $\Gamma$, and $\Pt{}$ denote $\frac d{d\ts}$. 
The geodesic curvature vector $\vect{k_g}$ is the image of the orthogonal projection of $\Ppt{\sigma}$ to $T_{\sigma}\La^4=(\spa{\sigma})^{\perp}$. 
As $\langle\sigma, \sigma\rangle=\langle\Pt{\sigma}, \Pt{\sigma}\rangle=1$, and moreover, $\langle\Ppt{\sigma}, \Ppt{\sigma}\rangle=1$ (\cite{La-So},\cite{Yu}), we have $\langle{\sigma}, \Ppt{\sigma}\rangle=-1$, and therefore, the geodesic curvature vector is given by $\vect{k_g} = \Ppt{\sigma}+\sigma$. 
It is lightlike, namely, $\Gamma$ is a {\em drill} (\cite{La-So}). 

\begin{figure}[ht] 
\begin{center}
\includegraphics[width=75mm]{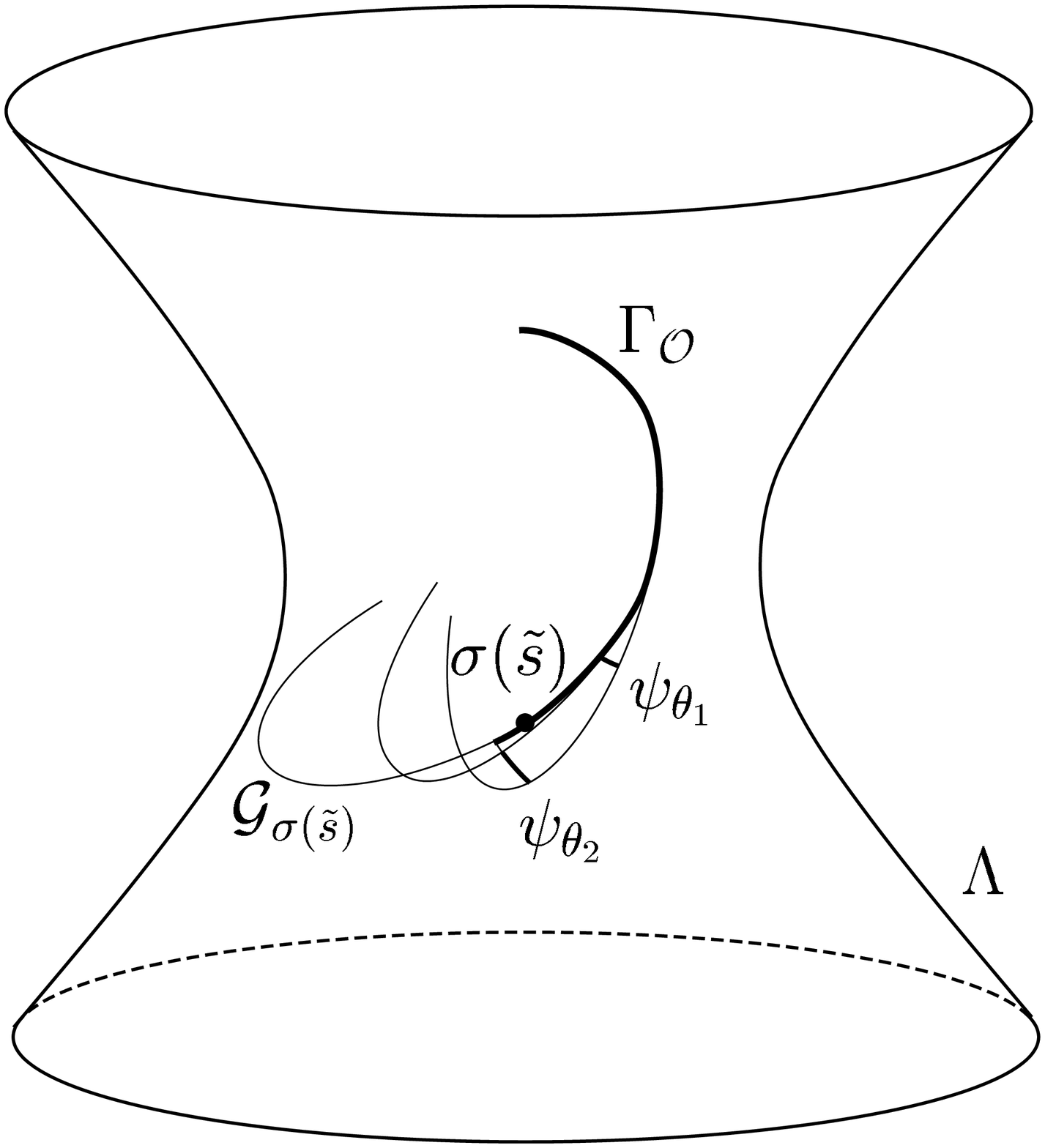}
\caption{$\Gamma_{\mathcal O}$ is a curve of osculating spheres $\Gamma_{\mathcal O}$, $\mathcal{G}_{\sigma(\ts)}$ is a geodesic circle in $\La$ that is tangent to $\Gamma_{\mathcal O}$ at $\sigma(\ts)$. 
Lightlike curves made of spheres containing the osculating circles are orthogonal to $\mathcal{G}_{\sigma(\ts)}$.}
\label{light_curves}
\end{center}
\end{figure}

Let $\mathcal{G}_{\sigma(\ts)}$ be the geodesic circle of $\La^4$ which is tangent to $\Ga$ at $\sigma(\ts)$ (Figure \ref{light_curves}). 
Then it is given by 
$\mathcal{G}_{\sigma(\ts)}=\La^4\cap\span{\sigma(\ts), \Pt{\sigma}(\ts)}.$ 
A point on it can be expressed as $\nu(\ts, \theta)=\cos\theta \sigma(\ts)+ \sin\theta \Pt{\sigma}(\ts)$ for some $\theta$. 
It corresponds to a sphere which contains the osculating circle to $C$ at $c(\ts)$. 

Let $V(C)=\bigcup_{\ts}\mathcal{G}_{\sigma(\ts)}$ be a surface in $\Lambda^4$ which is the union of the geodesic circles tangent to $\Gamma$:
\[V(C)=\{\nu(\ts, \theta)=\cos\theta \sigma(\ts)+ \sin\theta \Pt{\sigma}(\ts)\,|\,\ts\in I, 0\le\theta\le2\pi\}.\]
The tangent space of $V(C)$ at $\nu=\nu(\ts, \theta)$ is given by 
$$T_{\nu(\ts, \theta)}V(C)=T_{\nu}\Lambda^4 \cap \span{\sigma(\ts), \Pt{\sigma}(\ts), \Ppt{\sigma}(\ts)}= T_{\nu}\Lambda^4 \cap \span{\sigma(\ts), \Pt{\sigma}(\ts), \vect{k_g} (\ts)}.$$ 
As 
\[T_{\nu(\ts, \theta)}\mathcal{G}_{\sigma(\ts)}=\span{-\sin\theta\sigma(\ts)+\cos\theta\Pt{\sigma}(\ts)}\]
its orthogonal complement in $T_{\nu}V(C)$ is given by 
$\span{\Ppt{\sigma}+\sigma}=\span{\vect{k_g}(\ts)}$. 

It follows that the curves orthogonal to the ``foliation" of $V(C)$ by these geodesic circles $\{\mathcal{G}_{\sigma(\ts)}\}$ are lightlike. 
Generically they have cuspidal edges at points of $\Gamma$. 

Notice that the angle between two spheres  $\nu(s,\theta)$ and $\nu(s,\theta')$ is independant of the value of $s$. 

\begin{remark}\rm 
We have $d\ts=|T|c^{\ast}d\rho$, where $\rho$ is the conformal arc-length and $T$ is the {\em conformal torsion} given by 
\[T=\frac{2\kappa'\tau+\kappa^2\tau^3+\kappa\kappa'\tau'-\kappa\kappa''\tau}{(\kappa'^2+\kappa^2\tau^2)^{5/2}}\]
(see \cite{CSW}). 
\end{remark}

\begin{proposition}\label{prop_second_order}
Let $C$ be a vertex-free curve, $\overline{\Sigma}_0$ a sphere which has the second order contact with $C$, and $\bar{\sigma}_0$ a point in $\La$ corresponding to $\overline{\Si}_0$. 
Then there is a unique lightlike curve $\bar{\sigma}$ through $\bar{\sigma}_0$ consisting of the spheres $\overline{\Sigma}$  with the second order contact with $C$. 
\end{proposition}
The statement without ``{\sl uniqueness}'' was proved in Corollary 10 of \cite{La-So}. 

We remark that a sphere has the second order contact with a curve at a point $m$ if and only if it contains the osculating circle of the curve at $m$. 
\begin{proof} Every sphere having the second order contact at a point $m(\tilde{s})$ in $C$ can be written by 
\[
\cos \theta\, \sigma(\tilde{s})-\sin \theta\, \Pt{\sigma}(\tilde{s})
\]
for some $\theta\in[0,2\pi)$, where $\Pt{}$ denotes $d/d\tilde{s}$. 
Let $\eta(\tilde{s})$ be a curve in $\La$ given by 
\[
\eta(\tilde{s})=\cos u(\tilde{s})\, \sigma(\tilde{s})-\sin u(\tilde{s})\, \Pt{\sigma}(\tilde{s}), 
\]
where $u(\tilde{s})$ is a function. 
Then we have 
\[
\Pt{\eta}=\cos u\,\big(1-\Pt{u}\big)\,\Pt{\sigma}-\sin u\,\big(\Pt{u}\sigma+\Ppt{\sigma}\big).
\]
As $\langle\sigma, \sigma\rangle=\langle\Pt{\sigma}, \Pt{\sigma}\rangle=1$ we have $\langle\Pt{\sigma}, \sigma\rangle=\langle\Pt{\sigma}, \Ppt{\sigma}\rangle=0$. 
Furthermore we have $\langle\Ppt{\sigma}, \Ppt{\sigma}\rangle=1$ (\cite{La-So},\cite{Yu}). 
Therefore we have 
\[
\langle \Pt{\eta}, \Pt{\eta}\rangle=\big(1-\Pt{u}\big)^2, 
\]
which implies that $\eta(\tilde{s})$ is lightlike if and only if $u(\tilde{s})=\tilde{s}+\theta$ for some constant $\theta$. 

If $\bar{\sigma}_0$ can be expressed as 
\[
\bar{\sigma}_0=\cos u_0\, \sigma(\tilde{s}_0)-\sin u_0\, \Pt{\sigma}(\tilde{s}_0)
\]
then $\bar{\sigma}$ is uniquely determined by 
\[
\bar{\sigma}(\tilde{s})=\cos \,(\tilde{s}+u_0-\tilde{s}_0)\, \sigma(\tilde{s})-\sin \,(\tilde{s}+u_0-\tilde{s}_0)\, \Pt{\sigma}(\tilde{s}), 
\]
which has the following lightlike tangent vector
\begin{equation}
\Pt{\bar{\sigma}}(\tilde{s})=-\sin(\tilde{s}+u_0-\tilde{s}_0)\big(\sigma(\tilde{s})+\Ppt{\sigma}(\tilde{s})\big).
\end{equation}
\end{proof}

\smallskip
\begin{remark}
Let 
\[\psi(\theta,\tilde{s})=\cos (\tilde{s}+\theta) \sigma(\tilde{s})-\sin (\tilde{s}+\theta)\Pt{\sigma}(\tilde{s}),\]
and $\psi_{\theta}$ be a lightlike curve $\{\psi(\theta,\,\overset{\mbox{\Huge .}}{}\,)\}$ in $\La^4$. 
Using one of the curves $\psi_{\theta}$ we can construct a surface $M_{\theta}$ containing $C$ such that $C$ is a line of principal curvature of $M_{\theta}$. 
Then the spheres of $\psi_{\theta}$ are the osculating spheres to $M_{\theta}$ at points of $C$. 
Any surface tangent to $M_{\theta}$ along $C$ has the same property.

When $C$ is contained in $\RR^3$, a developable $M_{\theta}$ is obtained as the envelopes of planes $P_{\theta}(\tilde{s}+h)$ satisfying:

\begin{itemize}
\item[-] the angle between $P_{\theta}(\tilde{s})$ and the osculating plane to $C$ at $c(\tilde{s})$ is $\theta$, 
\item[-] the derivative, with respect to an arc-lenth parameter $s$ on the curve $C$, of the angle between $P_{\theta}(\tilde{s}+h)$ and the osculating plane to $C$ at $c(\tilde{s}+h)$ is $-\tau(\tilde{s}+h)$, the opposite of the torsion of the curve $C$ at $c(\tilde{s}+h)$. 
\end{itemize}
\end{remark}
\bigskip
Note that $\psi_{\theta}$ is a curve in the surface $V(C)$ which intersects geodesic circles $\mathcal{G}_{\sigma(\ts)}$ orthogonally. 

\begin{corollary}\label{cor_average_lightlike}
The integral of the pull-back of the $\frac12$ dimensional length element {\rm (Definition \ref{def_1/2-length_element})} of the lightlike curves $\psi_{\theta}$, $\psi_{\theta}{}^{\ast}d\rho_{L^{\frac12}(\psi_{\theta})}$, with $\theta$ moving from $0$ to $2\pi$ is proportional to 
\begin{equation}\label{f_L(sigma''')-1}
\sqrt[4]{\big|\LL\big( \Pt{\sigma} +\Pppt{\sigma}\big)\big|}\,d\ts
=\sqrt[4]{\big|\LL\big(\Pppt{\sigma}\big)-1\big|}\,d\ts.
\end{equation}
\end{corollary}
\begin{proof}
Computing the second derivative of $\psi_{\theta}$, we get 
$$\Ppt{\psi}_{\theta}= -\cos(\ts+\theta)\left(\sigma(\ts)+\Ppt{\sigma}(\ts)\right)-\sin(\ts+\theta)\left(\Pt{\sigma}(\ts)+\Pppt{\sigma}(\ts)\right).$$ 
As $ \sigma(\ts)+\Ppt{\sigma}(\ts)$ is lightlike and orthogonal to its derivative $\Pt{\sigma}(\ts)+\Pppt{\sigma}(\ts)$, we have 
\[\sqrt[4]{\frac{\big|\LL\big(\Ppt{\psi}_{\theta}\big)\big|}{12}}
=\sqrt{|\sin(\ts+\theta)|} \,
\sqrt[4]{\frac{\big|\LL\big(\Pt{\sigma}(\ts)+\Pppt{\sigma}(\ts)\big)\big|}{12}}\,, \]
which implies 
\[\begin{array}{rcl}
\displaystyle \int_0^{2\pi}\psi_{\theta}{}^{\ast}d\rho_{L^{\frac12}(\psi_{\theta})}\,d\theta
&=& \displaystyle  \int_0^{2\pi}\sqrt[4]{\frac{\big|\LL\big(\Ppt{\psi}_{\theta}\big)\big|}{12}}\,d\ts\,d\theta \\[2mm]
&=& \displaystyle \left(\int_0^{2\pi}\sqrt{|\sin\theta|}\,d\theta\right)\sqrt[4]{\frac{\big|\LL\big(\Pt{\sigma}(\ts)+\Pppt{\sigma}(\ts)\big)\big|}{12}}\,d\ts\,.\end{array}\]
It gives the left hand side of (\ref{f_L(sigma''')-1}). 

On the other hand, the equation (\ref{f_L(sigma''')-1}) follows from Table \ref{table_sigma_i-sigma_j} of $\langle {\sigma}^{(i)}, {\sigma}^{(j)} \rangle$ $(0\le i,j\le 3)$, which can be obtained from $\big\langle\Ppt{\sigma}, \Ppt{\sigma}\big\rangle=1$ (\cite{La-So},\cite{Yu}). 

\begin{table}[h]
\begin{center}
\par\noindent
\begin{tabular}{|c|c|c|c|c|}
\hline
& $\!\!\phantom{\stackrel{{a}}{0}}\!\! \sigma\!\!$  &  $\!\!\Pt{\sigma}\!\!$  &  $\!\!\Ppt{\sigma}\!\!$  & $\!\!\Pppt{\sigma}\!\!$ \\[1mm] 
\hline
 $\!\!\phantom{\stackrel{{a}}{0}}\!\!\sigma$ & $1$ & $0$ & $-1$ & $0$ \\[1mm] 
\hline
 $\!\!\phantom{\stackrel{{a}}{0}}\!\!\Pt{\sigma}$ &  & $1$ & $0$ & $-1$  \\[1mm] 
\hline
 $\!\!\phantom{\stackrel{{a}}{0}}\!\!\Ppt{\sigma}$ &  &  & $1$ & $0$ \\[1mm] 
\hline
 $\!\!\phantom{\stackrel{{a}}{0}}\!\!\Pppt{\sigma}$ &  &  &  & $\LL\big(\Pppt{\sigma}\big)$ \\[1mm]
\hline 
\end{tabular}
\end{center}
\caption{A table of $\langle {\sigma}^{(i)}, {\sigma}^{(j)} \rangle$ }
\label{table_sigma_i-sigma_j}
\end{table}
\end{proof}

\begin{theorem}\label{average_light}
Let $C$ be a vertex-free curve. 
Let $\Ga=\{\sigma(\ts)\}$ be a curve in $\La^4$ which is the set of osculating spheres of $C$, where $\ts$ is the arc-length of $\Ga$. 
Then the conformal arc-length of the curve $C$ is, up to the multiplication by a universal constant, equal to the average with respect to $\theta$ of the $L^{\frac12}$-measures of the lightlike curves in 
\[V(C)=\{\nu(\ts, \theta)=\cos\theta \sigma(\ts)+ \sin\theta \Pt{\sigma}(\ts)\,|\,\ts\in I, 0\le\theta\le2\pi\}\]
through $\cos\theta \, \sigma(\ts_0) + \sin\theta \, \Pt{\sigma}(\ts_0)$. 
\end{theorem}
\begin{proof}
It is enough to show that the pull-back of the conformal arc-length element is proportional to the average of that of the $\frac12$ dimensional length element. 

Let $\Ga_{\mathcal O} \subset \Lambda^4$ be the curve of osculating spheres to the curve $C$, and $\sigma(\tilde{s})$ the osculating sphere at the point $m(\tilde{s})$ to the curve $C$, where $\tilde{s}$ is an arc-lenth parameter on the curve $\Ga_{\mathcal O} \subset \Lambda^4$. 
We remark that $\Ga_{\mathcal O}$ is spacelike if $C$ is vertex-free. 

Let $S$ be the sphere orthogonal to $C$ at $m(s)$ which also contains the point $m(s+h)$. 
 We will now follow the intersection with  $S$ of the spheres of the different lightlike curves $\psi_{\theta}$, where 
$\theta$ is the angle of the initial sphere of the family $\psi_{\theta}(t_0)$ with a chosen sphere of the second order contact with $C$ at the point $m_0 =c(s_0)$. 

Consider two points $m(s)$ and $m(s+h)$ on $C$ (here we use the arc-length parameter on the curve $C\subset \ss^3$). 
The set of the spheres $\{\psi(\theta, {s})\, | \theta \in \ss^1\}$ and $\{\psi(\theta, {s}+h)\, | \theta \in \ss^1\}$ are two pencils consisting of the spheres which contain ${\cal O}_{m(s)}$ and ${\cal O}_{m(s+h)}$ respectively. They intersect the sphere $S$ in two pencils of circles with base points $\{m(s),y(s+h)\}$ and $\{y(s),m(s+h)\}$  (see Figure \ref{trace_on_S}).

As the sphere $S$ is orthogonal to the osculating circle to the curve $C$ at $m(s)$, which is the base circle of the pencil $\{\psi(\theta,s)\, | \theta \in \ss^1\}$, the angle $\theta$ is also the angle parameter of the pencil of circles $\{\psi(\theta,s)\cap S\, | \theta \in \ss^1\}$. The angle of the circles $S\cap \psi(\theta,s+h)$ and $S\cap \psi(\theta',s+h)$ is only close to $\theta-\theta'$ as we know that the angle of $S$ and ${\cal O}_{m(s+h)}$ is close to $\pi/2$.

\begin{figure}[ht] 
\begin{center}
\includegraphics[width=8cm]{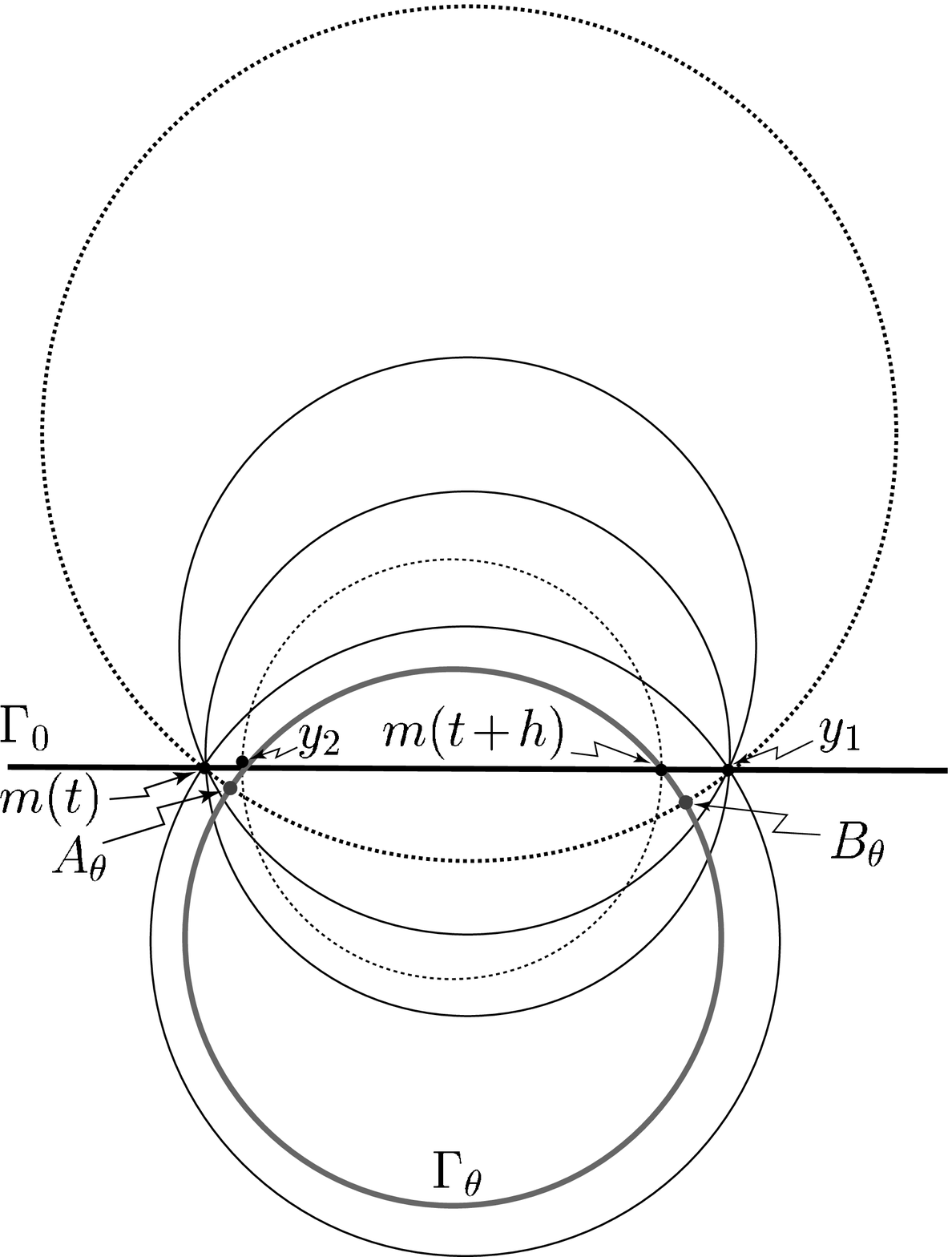}
\caption{Trace of the spheres $\psi_{\theta}({s})$ and $\psi_{\theta}({s}+h)$ on an almost orthogonal sphere $S$ (with a point of $S$ at infinity) \label{trace_light_curve_spheres}}
\end{center}
\end{figure}

Let ${\Gamma_0}$ denote the circle of ${S}$ containing the three points ${m}(s), y(s)$ an $m(s+h)$. As the cross ratio $cross({m(s)}, y(s+h), {m(s+h)}, y(s))$ is almost real (that is, the quotient of the imaginary part divided  by the real part is of order $o(h)$), the fourth point $y(s+h)$ is almost on ${\Gamma_0}$ (that is, after performing homothety which makes the radius of $S$ being equal to $1$ the distance between $y(s+h)$ and ${\Gamma_0}$ is of order $o(h)$). 
Therefore the cross ratio $cross({m(s)}, y(s+h), {m(s+h)}, y(s))$ is equivalent to the Lorentz distance between the two circles $S\cap \psi(\theta_1,s)$ and $S\cap \psi(\theta_1,s+h)$, where the value $\theta_1$ is such that $S\cap \psi_{\theta_1}(s)$ is orthogonal to ${\Gamma_0}$.  Let us now consider the intersection of the two circles $S\cap \psi(\theta,s)$ and $S\cap \psi(\theta,s+h)$ with the circle ${\Gamma}_{\theta}$ which is orthogonal to $S\cap \psi(\theta,s)$ (see Figure \ref{trace_light_curve_spheres}). 
The Lorentz distance between these two circles is equivalent to the cross ratio 
\[cross(m(s),A_{\theta}; B_{\theta},y_1)\simeq \cos^2(\theta -\theta_1)\, cross({m(s)}, {y_2}, {m(s+h)}, {y_1}). \]
Using Lemma \ref{cross_almost}, we see that the Lorentz distance between the spheres $\psi_{\theta}(s)$ and $\psi_{\theta}(s+h)$ is also equivalent to 
\[\cos(\theta -\theta_1)\, cross({m(s)}, y(s+h), {m(s+h)}, y(s)).\]

The Lorentz distance between two close disjoint spheres is equivalent to their Lorentz distance. We can integrate the contribution of the segments $\psi_{\theta}(t), \psi_{\theta}(s+h)$ to the $L^{1/2}$-lengthes of the curves $\psi_{\theta} \subset \Lambda$. It is 
 $$
\displaystyle  \left( \int_{\theta \in \ss^1} |\cos(\theta -\theta_1)|^{1/2}d\theta \, \right) \, \sqrt{2}\,\sqrt[4]{|cross({m(s)}, y(s+h), {m(s+h)}, y(s))|}\,.
$$
\noindent This provides the constant in the statement of Theorem \ref{average_light}.
\end{proof}

\medskip
\begin{remark}\rm 
There is an alternative proof of Theorem \ref{average_light} using (\ref{f_L(sigma''')-1}) and an anti-isometry between two Grassmann manifolds. 

Recall that the set $\mathcal S(1,3)$ of oriented circles in $\SS^3$ can be identified with the Grassmann manifold $\widetilde{\textsl{Gr}}_-(3;\mathbb{R}^{5}_1)$ of oriented $3$-dimensional timelike subspaces of $\mathbb{R}^{5}_1$. 
The orthogonal complement of a timelike $3$-space is a spacelike $2$-space. 
Therefore, there is a bijection between $\widetilde{\textsl{Gr}}_-(3;\mathbb{R}^{5}_1)\subset\RR^{10}_4$ and the Grassmann manifold $\widetilde{\textsl{Gr}}_+(2;\mathbb{R}^{5}_1)\subset\RR^{10}_6$ of oriented $2$-dimensional spacelike subspaces of $\mathbb{R}^{5}_1$ that can be obtained by assigning the orthogonal complement. 
This bijection is a restriction of an anti-isometry $F$ (here $F$ being an anti-isometry means $\langle F(u), F(v)\rangle=-\langle u, v\rangle$ for any $u, v$) between $\RR^{10}_4$ and $\RR^{10}_6$ that exchanges spacelike and timelike subspaces (\cite{La-OH2}). 

Through this anti-isometry, the osculating circle $\gamma$ corresponds to $\sigma\w\Pt{\sigma}$ in $\widetilde{\textsl{Gr}}_+(2;\mathbb{R}^{5}_1)$, and therefore, $\Ppt{\gamma}$ corresponds to $\sigma\w\Pppt{\sigma}+\Pt{\sigma}\w\Ppt{\sigma}$. 
Then, by Table \ref{table_sigma_i-sigma_j} we have $\LL\big(\Ppt{\gamma}\big)=\LL\big(\Pppt{\sigma}\big)-1$, which implies that the pull-back of the $\frac12$ dimensional length element of $\gamma$, and hence the pull-back of the conformal arc-length element, is given by $\sqrt[4]{\big|\LL\big(\Pppt{\sigma}\big)-1\big|}\,d\ts$. 
Now Corollary \ref{cor_average_lightlike} implies that the integral of the pull-back of the $\frac12$ dimensional length element of the lightlike curves $\psi_{\theta}$, $\psi_{\theta}{}^{\ast}d\rho_{L^{\frac12}(\psi_{\theta})}$, with $\theta$ moving from $0$ to $2\pi$ is proportional to the pull-back of the conformal arc-length element, which completes the proof. 
\end{remark}

\bigskip
The reader is referred to \cite{O'N} for pseudo-Riemannian metrics, 
\cite{Ak-Go}, \cite{Fi}, \cite{HJ}, \cite{MRS}, \cite{Ro-Sa}, and \cite{Su1, Su2} for further details in conformal differential geometry, and \cite{Yu} for the construction of a conformally invariant moving frame along a curve in a spherical model in $\RR^5_1$. 

\bigskip \noindent
Institut de Math\'ematiques de Bourgogne, Universit\'e de Bourgogne, \\
CNRS-UMR 5584, U.F.R. Sciences et Techniques, 9, avenue Alain Savary - B.P. 47870, 21078 DIJON Cedex, FRANCE.\\
E-mail: langevin@u-bourgogne.fr

\bigskip \noindent
Department of Mathematics, Tokyo Metropolitan University, \\
1-1 Minami-Ohsawa, Hachiouji-Shi, Tokyo 192-0397, JAPAN. \\
E-mail: ohara@tmu.ac.jp
\end{document}